\documentclass[a4paper]{amsart}
\newif\ifams
\amstrue
\usepackage{amsmath}
\usepackage{mathtools}%
\usepackage{booktabs}
\usepackage{amssymb}
\usepackage{enumerate}
\usepackage{bm}%
\usepackage{natbib}%
\usepackage{macros}
\makeatletter
\DeclareRobustCommand*\cal{\@fontswitch\relax\mathcal}
\makeatother
\newenvironment{theorem}{\begin{Theorem}}{\end{Theorem}}
\newenvironment{definition}{\begin{Definition}}{\end{Definition}}
\newenvironment{lemma}{\begin{Lemma}}{\end{Lemma}}

\newenvironment{remark}{\begin{Remark}}{\end{Remark}}

\newenvironment{corollary}{\begin{Corollary}}{\end{Corollary}}

\renewcommand{\qedsymbol}{}
\def\keywords{\smallskip\noindent\textsc{Keywords. }}
\usepackage{tikz}
\usetikzlibrary{calc}
\usepackage{url,hyperref}
\usepackage[disable,textsize=scriptsize]{todonotes}
\begin{document}
\renewcommand{\h}{\bm h}
\renewcommand{\w}{\bm w}
\renewcommand{\o}{\bm o}
\title[On Ordinal Invariants in Well Quasi Orders]{On Ordinal Invariants in
  Well Quasi Orders\ifams\\\fi and Finite Antichain Orders}

\thanks{Mirna D{\v z}amonja thanks the Leverhulme Trust for a Research
  Fellowship for the academic year 2014/2015, the Simons Foundation
  for a Visiting Fellowship in the Autumn of 2015 and the Isaac Newton
  Institute in Cambridge for their support during the HIF programme in
  the Autumn of 2015, supported by EPSRC Grant Number
  EP/K032208/1.  She most gratefully acknowledges the support of the
  {\em Institut d'histoire et de philosophie des sciences et des
    techniques (IHPST)} at the University Paris 1, where she is as an
  Associate Member.  The three authors thank the London Mathematical
  Society for their support through a Scheme 3 Grant.}
\author[M.~D{\v z}amonja]{Mirna D{\v z}amonja}
\address{University of East Anglia, Norwich, UK}
\email{m.dzamonja@uea.ac.uk}
\author[S.~Schmitz]{Sylvain Schmitz}
\address{LSV, ENS Paris-Saclay \& CNRS, Universit\'e Paris -
  Saclay, France}
\email{\{schmitz,phs\}@lsv.ens-cachan.fr}
\author[Ph.~Schnoebelen]{Philippe Schnoebelen}

\begin{abstract}
We investigate the ordinal invariants height, length, and width of
well quasi orders (WQO), with particular emphasis on width, an
invariant of interest for the larger class of orders with finite
antichain condition (FAC).  We show that the width in the class of FAC
orders is completely determined by the width in the class of WQOs, in
the sense that if we know how to calculate the width of any WQO then
we have a procedure to calculate the width of any given FAC order.  We
show how the width of WQO orders obtained via some classical
constructions can sometimes be computed in a compositional way.  In
particular this allows proving that every ordinal can be
obtained as the width of some WQO poset.
One of the difficult questions is to give a complete formula for the
width of Cartesian products of WQOs.  Even the width of the product of
two ordinals is only known through a complex recursive formula.
Although we have not given a complete answer to this question we have
advanced the state of knowledge by considering some more complex
special cases and in particular by calculating the width of certain
products containing three factors.
In the course of writing the paper we have discovered that some of the
relevant literature was written on cross-purposes and some of the
notions re-discovered several times.  Therefore we also use the
occasion to give a unified presentation of the known results.

\keywords wqo, width of wqo, ordinal invariants
\end{abstract}
\maketitle

\section{Introduction}
\label{intro}

In the finite case, a partial order---also called a
\emph{poset}---$(P,{\leq})$ has natural cardinal invariants: a
\emph{width}, which is the cardinal of its maximal antichains, and a
\emph{height}, which is the cardinal of its maximal chains.  The width
and height are notably the subject of the theorems of \citet{Dilworth}
and \citet{Mirsky} respectively; see \citet{west82} for a survey of
these \emph{extremal} problems.  In the infinite case, cardinal
invariants are however less informative---especially for countable
posets---,\linebreak while the theorems of \citeauthor{Dilworth} and
\citeauthor{Mirsky} are well-known to fail%
~\citep{Peles,Schmerl}.

When the poset at hand enjoys additional conditions, the corresponding
\emph{ordinal invariants} offer a richer theory, as studied for instance by
\citet{kriz90b}.  Namely, if $(P,{\leq})$ has the \emph{the finite
  antichain condition} (FAC), meaning that its antichains are finite,
then the tree
\begin{align*}
  \Inc(P)&\eqdef\bigl\{\langle x_0,x_1,\dots,x_n\rangle \in P^{<\omega}
  ~:~
  0\leq n<\omega \land \forall 0\leq i<j\leq n,\,x_i\mathbin\bot
  x_j\bigr\}
  \intertext{%
    of all non-empty (finite) sequences of pairwise
    \underline{inc}omparable elements of $P$ ordered by initial
    segments has no infinite branches.  Note that the tree
    $(\Inc(P),{\initial})$ does not necessarily have a single root and
    that the empty sequence is excluded (the latter is a matter of
    aesthetics, but it does make various arguments run more smoothly
    by not having to consider the case of the empty sequence
    separately).  Therefore, $\Inc(P)$ has a rank, which is the
    smallest ordinal $\gamma$ such that there is a function
    $f:\,\Inc(P)\to\gamma$ with $s\mathrel\initial t\implies f(s)> f(t)$ for
    all $s,t\in \Inc(P)$.  This ordinal is called the \emph{width} of
    $P$ and in this paper we denote it by $\w(P)$---it was
    denoted by $\mathrm{wd}(P)$ by \citet{kriz90b}.
    \newline\hspace*{\parindent}
    Similarly, if $(P,{\leq})$ is \emph{well-founded} (WF), also
    called {Artinian},
    meaning that its descending sequences are finite, then the tree}
  \Dec(P)&\eqdef\bigl\{\langle x_0,x_1,\dots,x_n\rangle \in P^{<\omega}
  ~:~
  0\leq n<\omega \land \forall 0\leq i<j\leq n,\,x_i> x_j\bigr\}
  \intertext{%
    of non-empty strictly descending sequences has an ordinal rank, which we
    denote by $\h(P)$ (\citeauthor{kriz90b} denote it by
    $\mathrm{ht}(P)$) and call the \emph{height} of $P$.
    \newline\hspace*{\parindent}
    Finally, if $(P,{\leq})$ is both well-founded and FAC, i.e., is a
    \emph{well partial order} (WPO), then the tree}
  \Bad(P)&\eqdef\bigl\{\langle x_0,x_1,\dots,x_n\rangle \in P^{<\omega}
  ~:~
  0\leq n<\omega \land \forall 0\leq i<j\leq n,\,x_i\not\leq
  x_j\bigr\}
\end{align*}
of non-empty \emph{bad sequences} of $P$ has an ordinal rank, which we
denote by $\o(P)$ and call the \emph{maximal order type} of $P$ after
\citet{deJonghParikh} and \citet{schmidt79} (\citeauthor{kriz90b}
denote it by $c(P)$, \citeauthor{BlGu} call it the \emph{stature} of
$P$).  In the finite case, this invariant is simply the cardinal of
the poset.

Quite some work has already been devoted to heights and maximal order
types, and to their computation.  Widths are however not that
well-understood: as \citet[Rem.~4.14]{kriz90b} point out, they do not
enjoy nice characterisations like heights and maximal order types do,
and the range of available results and techniques on width
computations is currently very limited.

\medskip
Our purpose in this paper is to explore to what extent we can find
such a characterisation, and provide formul\ae\ for the behaviour of
the width function under various classically defined operations with
partial orders.  Regarding the first point, we first show in
Sec.~\ref{sec-characterizations} that the width coincides with the
\emph{antichain rank} defined by \citet{AbBo}, which is the height of
the chains of antichains; however, unlike the height and maximal order
type of WPOs, the width might not be attained
(Rem.~\ref{rk-max-width}).  Regarding the second point, we first show
in Sec.~\ref{wqoversusall} that computing widths in the class of FAC
orders reduces to computing widths in the class of WPOs.  We recall
several techniques for computing ordinal invariants, and apply them in
Sec.~\ref{sec-computing-w} to obtain closed formul\ae\ for the width
of sums of posets, and for the finite multisets, finite sequences, and
tree extensions of WPOs.  One of the main questions is to give a
complete formula for the width of the Cartesian products of WPOs.
Even the width of the product of two ordinals is only known through a
complex recursive formula (due to Abraham, see Sec.~\ref{finiteproducts}) and we only
have partial answers to the general question.

The three ordinal invariants appear in different streams of the
literature, often unaware of the results appearing in one another, and
using different definitions and notations.  Another motivation of this
paper is then to provide a unified presentation of the state of the
knowledge on the subject, and we also recall the corresponding results
for heights and maximal order types as we progress through the paper.

\section{Background and Basic Results}

\subsection{Posets and Quasi-Orders}
We consider posets and, more generally, quasi-orders (QO).  When
$(Q,{\leq_Q})$ is a QO, we write $x<_Qy$ when $x\leq_Q y$ and
$y\not\leq_Q x$.  We write $x\perp_Q y$ when $x\not\leq_Q y$ and
$y\not\leq_Q x$, and say that $a$ and $b$ are \emph{incomparable}.  We
write $x\equiv_Q y$ when $x\leq_Q y\land y\leq_Q x$: this is an
equivalence and the quotient $(Q,{\leq_Q})/\equiv_Q$ is a poset that,
as far as ordinal invariants are concerned, is indistinguishable
from~$Q$.  Therefore we restrict our attention to posets for technical
reasons but without any loss of generality.  Note that some
constructions on posets (e.g., taking powersets) yield quasi-orders
that are not posets.  A QO $Q$ is \emph{total} if for all $x,y$ in
$Q$, $x\le_Q y$ or $x\ge_Q y$; a total poset is also called a
\emph{chain}.

When a QO does not have infinite antichains, we say that it satisfies
the \emph{Finite Antichain Condition}, or simply that it is FAC.  A QO
that does not have any infinite (strictly) decreasing sequence is said
to be well-founded (or WF).  A \emph{well-quasi order} (or WQO) is a
QO that is both WF and FAC: it is well-known that a QO is WQO if and
only if it does not have any infinite bad
sequence~\cite{Kruskal,Milner}, where a sequence $\langle
x_0,x_1,x_2,\ldots\rangle$ is \emph{good} if $x_i\leq x_j$ for some
positions $i<j$, and is \emph{bad} otherwise.

For a QO $(Q, {\le})$ we define the \emph{reverse} QO $Q^\ast$ as $(Q,
{\ge})$, that is to say, $x\le_{Q^\ast} y$ if and only if $x\ge_Q y$.
An \emph{augmentation} of $(Q, {\le})$ is a QO $(Q, {\le'})$ such that
$x\le y\implies x\le' y$, i.e., $\le$ is a subset of $\le'$.  A
\emph{substructure} of a QO $(Q, {\le})$ is a QO $(Q', {\le'})$ such
that $Q'\subseteq Q$ and ${\le'}\:\subseteq\: {\le}$. In this case, we
write $Q'\le Q$.

\subsection{Rankings and Well-Founded Trees}\label{ssec-rank}
Recall that for every WF poset $P$ there exist ordinals $\gamma$ and
order preserving functions $f{:}\,P\to \gamma$, that is, such that
$x<_P y\implies f(x)< f(y)$ for all $x,y\in P$.  The smallest such
ordinal $\gamma$ is called the \emph{rank} of $P$; one can obtain the
associated \emph{ranking function} $r{:}\,P\to \gamma$ by defining
inductively $r(x)=\sup\{r(y)+1:\,y<_P x\}$, and the rank turns out to
be equal to its height $\h(P)$ (see Sec.~\ref{ssec:residuals}).  When $P$
is total, i.e., is a chain, then its rank is also called its
\emph{order type}.

Traditionally, for a tree $(T,\le_T )$, one says that it is
well-founded if it \emph{does not have an infinite branch}, which with
the notation above amounts to saying that the reverse partial order
$(T,\ge_T)$ is well-founded.  This somewhat confusing notation,
implies that for rooted well-founded trees, the root(s) have the
largest rank, and the leaves have rank~$0$.  In our definitions of
ordinal invariants given in the introduction, we considered trees of
non-empty finite sequences, ordered by initial segments: if $s=\langle
x_0,x_1,\ldots,x_n\rangle$ and $t=\langle y_0,y_1,\ldots,y_m\rangle$,
we write $s\initialeq t$ and say that $s$ is an \emph{initial segment}
of $t$, when $n\leq m$ and $s=\langle y_0,\ldots,y_n\rangle$.
Equivalently, the associated strict ordering $s\initial t$ means that
$t$ can be obtained by appending some sequence $t'$ after $s$, denoted
$t=s\frown t'$.

We also make an easy but important observation regarding
substructures: When $P$ is embedded in $Q$ as an induced substructure,
then $\w(P)\le \w(Q)$, and similarly for $\o$ and $\h$.  Indeed, every
antichain (bad sequence, decreasing sequence, resp.) of $P$ is an
antichain (bad sequence, decreasing sequence, resp.) of $Q$, so the
ranks of the corresponding trees can only increase when going from $P$
to~$Q$.

\subsection{Residual Characterisation}
\label{ssec:residuals}

For a poset $(P,{\le})$, $x\in P$, and
$\ast\in\{{\bot},{<},{\not\ge}\}$, we define the
\emph{$\ast$-residual} of $P$ at $x$ as the induced poset defined by
\begin{equation}
  P_{\ast x}\eqdef \{y\in P~:~ y\mathrel\ast x\}\;.
\end{equation}
Since this is an induced substructure of $P$, $P_{\ast x}$ is FAC
(resp.\ WF, WPO) whenever $P$ is FAC (resp.\ WF, WPO).

The interest of $\bot$-residuals (resp.\ $<$-residuals,
$\not\ge$-residuals) is that they provide the range of choices for
continuing incomparable (resp.\ descending, bad) sequences once
element $x$ has been chosen as first element: the suffix of the
sequence should belong to $P_{\ast x}$, and we have recursively
reduced the problem to measuring the rank of the tree $\Inc(P_{\bot
  x})$ (resp.\ $\Dec(P_{<x})$, $\Bad(P_{\not\ge x})$).

The following lemma shows precisely how we can extract the rank from
such a recursive decomposition of the tree.
\begin{lemma}\label{bunching}
  \hfill\begin{enumerate}
\item Suppose that $\{T_i:\,i\in I\}$ is a family of
well-founded trees and let $T$ be their disjoint union. Then
$T$ is a well-founded tree and it has rank
$\rho(T)=\sup_{i\in I} \rho(T_i)$.
\item Let $T=t^\frown F$ denote a tree rooted at $t$ with $F =
  T\setminus t$ and suppose that $F$ is well-founded of rank
  $\rho(F)$. Then so is $T$, and $\rho(T) =\rho(F) +1$.
\end{enumerate}
\end{lemma}

\begin{proof}[\ifams Proof \fi of~1] It is clear that $T$ is well founded.
For each $i\in I$, let $f_i{:}\,T_i\to
\rho(T_i)$ be a function witnessing the rank of $T_i$. Then
$f\eqdef\bigcup_{i\in I}f_i$ is an order reversing
function from $T$ to $\gamma\eqdef\sup_{i\in I}
\rho(T_i)$, showing $\rho(T)\le \gamma$.

Conversely, if $f{:}\,T\to\rho(T)$ is a witness function for
the rank of $T$, its restriction to any $T_i$ is
order reversing, showing that $\rho(T_i)\leq\rho(T)$.
\end{proof}
\begin{proof}[\ifams Proof \fi of~2]\renewcommand{\qedsymbol}{$\eop_{\ref{bunching}}$}
  Clearly $T$ is
  well-founded. Let $\rho^\ast\eqdef\rho(F)+1=\bigl(\sup_{\alpha<\rho(F)}(\alpha+1)\bigr)
  +1$.  Consider the ranking function $r{:}F\to \rho(F)$, and let
  $f{:}T\to\rho^\ast$ be given by
\[
  f(s) \eqdef \begin{cases}
          r(s)& \textrm{if }s\in F\:,      \\
          \sup_{\alpha< \rho(F)}(\alpha+1)& \textrm{if }s=t.
         \end{cases}
\]
It is clear that $f$ is an order reversing function, witnessing 
$h(T)\leq\rho^\ast$. 
Suppose that $\beta<\rho^\ast$ and that $h{:}\,T\to\beta$
is an order reversing function.  In particular, $h(r) < f(r)$, so let
$\alpha< \rho(F)$ be such that $h(r)<\alpha+1$. Let $s\in F$ be such
that $f(s)=\alpha$. Hence $h(r)\le h(s)$, yet $r<_T s$, a
contradiction.
\ifams\relax\else\hfill$\eop_{\ref{bunching}}$\fi
\end{proof}

Lemma~\ref{bunching} yields the equations:
\ifams\begin{equation}\begin{aligned}\label{eq-w-decomp}
  \w(P)&=\sup_{x\in P}\{\w(P_{\bot x}) + 1\}\:,\\
  \h(P)&=\sup_{x\in P}\{\h(P_{{<}x}) + 1\}\:,\\
  \o(P)&=\sup_{x\in P}\{\o(P_{{\not\ge}x}) + 1\}\:,
\end{aligned}\end{equation}\else
\begin{align}\label{eq-w-decomp}
  \w(P)&=\sup_{x\in P}\{\w(P_{\bot x}) + 1\}\:,&
  \h(P)&=\sup_{x\in P}\{\h(P_{{<}x}) + 1\}\:, &
  \o(P)&=\sup_{x\in P}\{\o(P_{{\not\ge}x}) + 1\}\:,
\end{align}\fi
that hold for any FAC, WF, or WPO, poset $P$ respectively.
Note that it yields $\w(\emptyset)=\h(\emptyset)=\o(\emptyset)=0$.

Equation~\eqref{eq-w-decomp} is used very frequently
in the literature and provides for a method for computing ordinal
invariants recursively, which we call the \emph{method of residuals}.

Equation~\ref{eq-w-decomp} further shows that the function
$r(x)\eqdef\h(P_{<x})$ is the optimal ranking function of~$P$.  Thus
$\h(P)$ is the rank of~$P$, i.e.\ the minimal $\gamma$ such that there
exists a strict order-preserving $f{:}\,P\to\gamma$ (recall
Sec.~\ref{ssec-rank}).

\subsection{Games for WQO Invariants}\label{ssec:games}

One limitation of the method of residuals is that it tends to produce recursive rather than
closed formul\ae, see, e.g., \citet{SS-icalp11}.
Another proof technique adopts a game-theoretical point of view. This
is based on \cite[\S3]{BlGu}, which in turn can be seen as an
application of a classical game for the rank of trees to the specific
trees used for the ordinal invariants.  We shall use this technique to
obtain results about special products of more than two orders, see for
example Thm.~\ref{cor-productofsquares}.

The general setting is as follows. For a WQO $P$ and an ordinal
$\alpha$, the game $G_{P,\alpha}^*$ ---where $*$ is one of
$\h,\o,\w$--- is a two-player game where positions are pairs
$(\beta,S)$ of an ordinal and a sequence over $P$. We start in the
initial position $(\alpha,\langle\rangle)$. At each turn, and in
position $(\beta,S)$, Player 1 picks an ordinal $\beta'<\beta$ and
Player 2 answers by extending $S$ with an element $x$ from $P$. Player
2 is only allowed to pick $x$ so that the extended $S'=S\frown x$ is a
decreasing sequence (or a bad sequence, or an antichain) when $*=\h$
(resp.\ $*=\o$, or $*=\w$) and he loses the game if he cannot answer
Player 1's move. After Player 2's move, the new position is
$(\beta',S')$ and the game continues.  Player 2 wins when the position
has $\beta=0$ and hence Player 1 has no possible move.  The game
cannot run forever so one player has a winning strategy.  Applying
\cite[Prop.~23]{BlGu} we deduce that Player 2 wins in $G_{P,\alpha}^*$
iff $*(P)\geq\alpha$. As we are mostly interested in the invariant
$\w$, we shall adopt the notation $G_{P,\alpha}$ for~$G_{P,\alpha}^{\w}$.

\subsection{Cardinal Invariants}
We can connect the ordinal invariants with cardinal measures but this
does not lead to very fine bounds. Here are two examples of what can
be said.

\begin{lemma}\label{upperbound}
Suppose that $Q$ is a FAC quasi-order of cardinal
$\kappa\geq\aleph_0$. Then $\w(Q)<\kappa^+$, the cardinal successor of
$|Q|$.
\end{lemma}

\begin{proof}\renewcommand{\qedsymbol}{$\eop_{\ref{upperbound}}$}
The tree $\Inc(Q)$ has size equal to $\kappa$ and therefore its rank is
an ordinal $\gamma<\kappa^+$. 
\ifams\relax\else\hfill\qedsymbol\fi
\end{proof}

\begin{theorem}[Dushnik-Miller]\label{the-partitions}
Suppose that $P$ is a WPO of cardinal $\kappa\ge\aleph_0$.  Then
$\h(P)\ge \kappa$.
\end{theorem}

\begin{proof}\renewcommand{\qedsymbol}{$\eop_{\ref{the-partitions}}$}
This is an easy consequence of Thm. 5.25 in \cite{DushnikMiller}. By
the definition of $\h$, it suffices to show that $P$ has a chain of
size $\kappa$.  Define a colouring $c$ on the set $[P]^2$ of pairs of
$P$ by saying $c({x,y})\eqdef 0$ if $x$ is comparable to $y$ and
$c({x,y})\eqdef 1$ otherwise.  Then use the relation $\kappa\arrows
(\kappa, \aleph_0)^2$, meaning that $P$ has a chain of cardinal
$\kappa$ or an antichain of cardinal $\aleph_0$, which for
$\kappa=\aleph_0$ is the Ramsey Theorem, and for $\kappa>\aleph_0$ is
the Dushnik-Miller Theorem.  Since $P$ is FAC, we must have a chain of
order type at least $\kappa$.
\ifams\relax\else\hfill\qedsymbol\fi
\end{proof}

Such results are however of little help when the poset at hand is
countable, because they only tell us that the invariants are countable
infinite, as expected.  This justifies the use of ordinal invariants
rather than cardinal ones.

\subsection{WPOs as a Basis for FAC Posets}\label{wqoversusall}

A \emph{lexicographic sum} of posets in some family $\{ P_i:\, i \in
Q\}$ of disjoint orders \emph{along} a poset $(Q,\leq_Q)$, denoted by
$\sum_{i\in Q}P_i$, is defined
as the order $\le$ on the disjoint union $P$ of $\{ P_i:\, i \in Q\}$
such that for all $x,y \in P$ we have $x\le y$ iff $x,y\in P_i$ for
some $i\in Q$ and $x\le_{P_i} y$, or $x\in P_i$ and $y\in P_j$ for some
$i,j \in Q$ satisfying $i <_Q j$.

The lexicographic sum of copies of $P$ along $Q$ is denoted by $P\cdot
Q$ and called the \emph{direct product} of $P$ and $Q$.  The
\emph{disjoint sum} of posets in $\{ P_i:\,i \in Q\}$ is defined as
the union of the orders $\le_{P_i}$: this is just a special case of a
lexicographic sum, where the sum is taken over an antichain~$Q$.  In
the case of two orders $P_1,P_2$, the lexicographic sum is denoted by
$P_1\sqcup P_2$.

As a consequence of Thm.~7.3 of \citet{abcdzt} (by taking the union
over all infinite cardinals~$\kappa$), one obtains the following
classification theorem.
\begin{theorem}[\citeauthor{abcdzt}]
\label{allFAC}
Let $\BB\!\PP$ be the class of posets which are either a WPO, the
reverse of a WPO, or a linear order.  Let $\PP$ be the closure of
$\BB\!\PP$ under lexicographic sums with index set in $\BB\!\PP$ and
augmentation.  Then $\PP$ is exactly the class of all FAC posets.
\end{theorem}

We will use the classification in Thm.~\ref{allFAC} to see that if we
know how to calculate $\w(P)$ for $P$ an arbitrary WPO, then we can
bound $\w(P)$ for any FAC poset $P$. This in fact follows from some
simple observations concerning the orders in the class $\BB\!\PP$.

\begin{lemma}\label{basic}
(1) If $P$ is total, then $\w(P)=1$. In general, if all the
antichains in a poset $P$ are of length $\le n$ for some $n<\omega$,
then $\w(P)\le n$, and $\w(P) = n$ in the case that there are
antichains of length $n$.

\noindent
(2) For any poset $P$, $\Inc(P)=\Inc(P^\ast)$ and hence in the case of
FAC posets we have $\w(P^\ast)=\w(P)$.

\noindent
(3) If $P'$ is an augmentation of a FAC poset $P$, then $\Inc(P')$ is
a subtree of $\Inc(P)$ and therefore $\w(P')\le \w(P)$.

\noindent
(4) Let $P$ be the lexicographic sum of posets $\{ P_i:\,i \in L\}$
along some linear order $L$. Then $\Inc(P)=\bigcup_{i \in L}
\Inc(P_i)$ and in the case of FAC posets we have $\w(P)=\sup_{i \in L}
\w(P_i)$.
\end{lemma}

\begin{proof}\renewcommand{\qedsymbol}{$\eop_{\ref{basic}}$}
(1) The only non-empty sequences of antichains in a linear order $P$
  are the singleton sequences. It is clear that the resulting tree
  $\Inc(P)$ has rank~$1$, by assigning the value~$0$ to any singleton
  sequence. The more general statement is proved in the same way,
  namely if all the antichains in a poset $P$ are of length $< n$ for
  some $n<\omega$ then it suffices to define $f{:}\, \Inc(P) \to n$ by
  letting $f(s)\eqdef n - |s|$.

\noindent
(2), (3) Obvious.

\noindent
(4) This is the same argument as in Thm.~\ref{theorem-lexsum}.(3).
\ifams\relax\else\hfill\qedsymbol\fi
\end{proof}

In conjunction with Thm.~\ref{allFAC}, we conclude that the problem of
bounding the width of any given FAC poset is reduced to knowing how to
calculate the width of WQO posets.  This is the consideration of the
second part of this article, starting with Sec.~\ref{sec-computing-w}.

\section{Characterisations of Ordinal Invariants}
\label{sec-characterizations}

We recall in this section the known characterisations of ordinal
invariants.  With the method of residuals we can follow \citet{kriz90b} and
show that the height and maximal order types of WPOs also correspond
to their maximal chain heights (Sec.~\ref{ssec:chain}) and maximal
linearisation heights (Sec.~\ref{ssec:lin}), relying on results of
\citet{Wolk} and \citet{deJonghParikh} to show that these maxima are
indeed attained.  In a similar spirit, the width of a FAC poset is
equal to its antichain rank (Sec.~\ref{AB}), an invariant studied by
\citet{AbBo}---but this time it is not necessarily attained.  Finally,
in Sec.~\ref{sec-links} we recall an inequality relating all three
invariants and shown by \citet{kriz90b}.

\subsection{Height and Maximal Chains}\label{ssec:chain}

Given a WF poset $P$, let $\?C(P)$ denote its set of non-empty chains.
Each chain $C$ from $\?C(P)$ is well-founded and has a rank
$\h(C)$; we denote the supremum of these ranks by
$\mathrm{rk}_\?C P\eqdef\sup_{C\in\?C(P)}\h(C)$.  As explained for
example by
\citet[Thm.~4.9]{kriz90b}, we have
\begin{equation}\label{eq-sup-chain}
  \rk_\?C P \le \h(P)
\end{equation}
and this can be shown, for instance, by induction on the height using
the method of residuals.  Indeed, \eqref{eq-sup-chain} holds when
$P=\emptyset$, and for the induction step
\begin{align*}
  \sup_{C\in\?C(P)}\h(C)
  &\eqby{\eqref{eq-w-decomp}}\sup_{C\in\?C(P)}(\sup_{x\in
    C}\{\h(C_{<x})+1\})
  \leq\sup_{x\in P}\{ (\!\sup_{C'\in\?C(P_{<x})}\!\h(C')) + 1\}
  \intertext{because $C_{<x}$ is a chain in
    $\?C(P_{<x})$, and then by induction hypothesis~\eqref{eq-sup-chain}}
  \sup_{C\in\?C(P)}\h(C)&\leq\sup_{x\in P}\{\h(P_{<x})+1\}
  \eqby{\eqref{eq-w-decomp}}\h(P)\;.
\end{align*}

\begin{remark}\label{rk-wolk}
  The inequality in~\eqref{eq-sup-chain} can be strict.  For instance,
  consider the forest $F$ defined by the disjoint union $\{C_n:
  n\in\+N\}$ along $(\+N,{=})$, where each $C_n$ is a chain of height
  $n$, and add a new top element $t$ yielding $P\eqdef
  t^\frown F$.  Then $P$ is
  WF (but not FAC and is thus not a WPO).
  Note that $\h(P)=\h(F)+1=\omega+1$.  However, every
  chain $C$ in $\?C(P)$ is included in
  $t^\frown C_n$ for some $n$ and has height bounded by $n+1$,
  while $\rk_\?C (P)=\omega<\h(P)$.\hfill$\eop_{\ref{rk-wolk}}$
\end{remark}

\Citet[Thm.~9]{Wolk} further shows that, when $P$ is a WPO, the
supremum is attained, i.e.\ there is a chain $C$ with rank
$\h(C)=\mathrm{rk}_\?C P$.  In such a case, \eqref{eq-sup-chain} can
be strengthened to
\begin{equation}\label{eq-max-chain}
  \max_{C\in\?C(P)}\h(C) = \rk_\?C P =  \h(P)
\end{equation}
as can be checked by well-founded induction with
\begin{align*}
  \h(P)&\eqby{\eqref{eq-w-decomp}}\sup_{x\in P}\{\h(P_{<x})+1\}
  \leq \sup_{x\in P}\{\h(C_x)+1\}\leq \sup_{x\in P}\h(C_x\cup\{x\})\leq \sup_{C\in\?C(P)}\h(C)
\end{align*}
where $C_x$ is a chain of $P_{<x}$ witnessing~\eqref{eq-max-chain} by
induction hypothesis, and $C_x\cup\{x\}$ is therefore a chain in
$\?C(P)$ of height $\h(C_x)+1$.

\begin{theorem}[%
 \citeauthor{Wolk}; \citeauthor{kriz90b}]\label{thm-equivalences-2} Let
 $P$ be a WPO.
  Then $\h(P)=\mathrm{rk}_\?C P=\max_{C\in\?C(P)}\h(C)$ is the maximal
  height of the non-empty chains of~$P$.
\end{theorem}

More generally, the WPO condition in Thm.~\ref{thm-equivalences-2} can
be relaxed using the following result proven in
\citep{pouzet79,schmidt81,milner81}.
\begin{theorem}[\citeauthor{pouzet79}; \citeauthor{schmidt81};
    \citeauthor{milner81}]\label{thm-max-chain}
  Let $P$ be a WF poset.  Then
  \begin{itemize}
    \item \emph{either} $\rk_\?C P=\max_{C\in\?C(P)}\h(C)$,
      i.e.\ there exist chains of maximal height,
    \item \emph{or} there exists an antichain $A$ of $P$ such that the
      set of heights $\{\h(P_{{<}x}):x\in A\}$ is infinite.
  \end{itemize}
\end{theorem}

\subsection{Maximal Order Types and Linearisations}\label{ssec:lin}

A \emph{linearisation} of a poset $(P,{\leq})$ is an augmentation
$L=(P,{\preceq})$ which is a total order: $x\leq y$ implies
$x\preceq y$.  We let $\?L(P)$ denote the set of linearisations of
$P$.  As stated by \cite{deJonghParikh}, a poset is a WPO if and only
if all its linearisations are well-founded.
\Citeauthor{deJonghParikh} furthermore considered the supremum
$\sup_{L\in\?L(P)}\h(L)$ of the order types of the linearisations of
$P$, and showed that this supremum was attained
\citep[Thm.~2.13]{deJonghParikh}; this is also the subject of of
\cite[Thm.~10]{BlGu}.

\begin{theorem}[%
  \citeauthor{deJonghParikh}; \citeauthor{kriz90b}]\label{thm-equivalences-1} Let $Q$ be a WQO.
  Then $\o(Q)=\max_{L\in\?L(Q)}\h(L)$ is the maximal
  height of the linearisations of~$Q$.
\end{theorem}

\subsection{Maximal Order Types and Height of Downwards-Closed Sets}

A subset $D$ of a WQO $(Q,{\leq})$ is \emph{downwards-closed} if, for
all $y$ in $D$ and $x\leq y$, $x$ also belongs to~$D$.  We let
$\?D(Q)$ denote the set of downwards-closed subsets of~$Q$.  For
instance, when $Q=\omega$, $\?D(\omega)$ is isomorphic to
$\omega+1$.

It is well-known that a quasi-order $Q$ is WQO if and only if it
satisfies the descending chain condition, meaning that
$(\?D(Q),{\subseteq})$ is well-founded.  Therefore $\?D(Q)$ has a rank
$\h(\?D(Q))$ when $Q$ is WQO.  As shown by \citet[Prop.~31]{BlGu},
this can be compared to the maximal order type of~$Q$.
\begin{theorem}[\citeauthor{BlGu}]
  Let $Q$ be a WQO.  Then $\o(Q)+1=\h(\?D(Q))$.
\end{theorem}

\subsection{Width and Antichain Rank}
\label{AB}

\Citet{AbBo} consider a structure similar to the tree $\Inc(P)$ for
FAC posets $P$, namely the poset $\?A(P)$ of all non-empty
antichains of $P$.  In the case of a FAC poset, the poset $(\?A(P),
{\supseteq})$ is well-founded.  Let us call its height the
\emph{antichain rank} of $P$ and denote it by $\rk_\?A
P\eqdef\h(\?A(P))$; this is the smallest ordinal $\gamma$ such that
there is a strict order-preserving function from $\?A(P)$
to~$\gamma$.

In fact the antichain rank and the width function we study have the
same values, as we now show.  Thus one can reason about the width
$\w(P)$ by looking at the tree $\Inc(P)$ or at
$(\?A(P),{\supseteq})$, a different structure.

\begin{theorem}\label{equal}
Let $P$ be a FAC poset. Then $\w(P)= \rk_\?A P$.
\end{theorem}
\begin{proof}\renewcommand{\qedsymbol}{$\eop_{\ref{equal}}$}
Let $\gamma=\rk_\?A P$ and let $r{:}\, \?A(P)\to \gamma$ be such
that $S\supsetneq T\implies r(S) < r(T)$ for all non-empty antichains
$S,T$.  Define $f{:}\,\Inc(P) \to \gamma$ by letting for $s$ non-empty
$f(s)\eqdef r(S)$, where $S$ is the set of elements of $s$.  This
function satisfies $s\initial t \implies f(s)> f(t)$ and hence
$\w(P)\le\rk_\?A P$.

Conversely, let $\gamma=\w(P)$ and $f{:}\,\Inc(P)\to\gamma$ be such
that $s\initial t \implies f(s)>f(t)$.  For a non-empty antichain
$S\in \?A(P)$, observe that there exist finitely many---precisely
$|S|!$--- sequences $s$ in $\Inc(P)$ with support set $S$.  Call this
set $\Lin(S)$ and define $r{:}\,\?A(P)\to\gamma$ by $r(S)\eqdef
\min_{s\in\Lin(S)}f(s)$.  Consider now an antichain $S$ with $r(S) =
f(s)$ for some $s\in\Lin(S)$, and an antichain $T$ with $T\supsetneq
S$: then there exists an extension $t$ of $s$ in $\Lin(T)$, which is
therefore such that $f(s)>f(t)$, and hence $r(S)=f(s)>f(t)\geq r(T)$.
Thus $\w(P)\geq\rk_\?A P$.
\ifams\relax\else\hfill\qedsymbol\fi
\end{proof}

\begin{remark}\label{rk-max-width}
  The width $\w(P)$ is in general not attained, i.e., there might not
  exist any chain of antichains of height $\w(P)$.  First note that
  even when $P$ is a WPO, $(\?A(P),{\supseteq})$ is in general not a
  WPO, hence Thm.~\ref{thm-equivalences-2} does not apply.  In fact,
  examples of FAC posets where the width is not attained abound.
  Consider indeed any FAC poset $P$ with $\w(P)\ge\omega$, and any
  non-empty chain $C$ in $\?C(\?A(P))$.  As $C$ is well-founded for
  $\supseteq$, it has a minimal element, which is an antichain
  $A\in\?A(P)$ such that, for all $A'\neq A$ in $C$,
  $A'\subsetneq A$.  Since $P$ is FAC, $A$ is finite, and $C$ is
  therefore finite as well: $\h(C)<\omega$.~\hfill$\eop_{\ref{rk-max-width}}$
\end{remark}

\subsection{Relationship Between Width, Height and Maximal Order Type}
\label{sec-links}

As we have seen in the previous discussion, $\w(P)=\h(\?A(P))$ the
antichain rank (where antichains are ordered by reverse inclusion).
\Citet[Thm.~4.13]{kriz90b} proved that there is another connection
between the ordinal functions discussed here and the width function.

The statement uses natural products of ordinals.  Recall for this that
the Cantor normal form (CNF) of an ordinal $\alpha$
\[
\alpha=\omega^{\alpha_0}\cdot m_0 + \cdots + \omega^{\alpha_\ell}\cdot
m_\ell
\]
is determined by a non-empty decreasing sequence $\alpha_0>\alpha_1
\cdots >\alpha_\ell\ge 0$ of ordinals and a sequence of natural numbers
$m_i> 0$. Cantor proved that every ordinal has a unique
representation in this form.  Two well-known operations can be defined
based on this representation: the \emph{natural or Hessenberg sum}
$\alpha\oplus\beta$ is defined by adding the coefficients of the
normal forms of $\alpha$ and $\beta$ as though these were polynomials
in $\omega$. The \emph{natural or Hessenberg product}
$\alpha\otimes\beta$ is obtained when the normal forms of $\alpha$ and
$\beta$ are viewed as polynomials in $\omega$ and multiplied accordingly.

\begin{theorem}[K\v r\'i\v z and Thomas]\label{thm-oandh}
 For any WQO $(Q,\leq)$ the following holds:
  \begin{gather}
    \label{eq-KT-ineq}
             \w(Q)\leq \o(Q) \leq \h(Q)\otimes \w(Q)\:.
  \end{gather}
\end{theorem}

For completeness, we give a detailed proof.
\begin{proof}\renewcommand{\qedsymbol}{$\eop_{\ref{thm-oandh}}$}
  For the first inequality, clearly any antichain in $Q$ can be
  linearised in an arbitrary way in a linearisation of $Q$. So $\w(Q)$
  is certainly bounded above by the length of the maximal such
  linearisation, which by Thm.~\ref{thm-equivalences-1} is exactly the
  value of~$\o(Q)$.

For the second inequality, let $\alpha=\w(Q)$ and let $g:\Inc(Q)\into
\alpha$ be a function witnessing that. Also, let $\beta=\h(Q)$ and let
$\rho:\,Q\into \beta$ be the rank function.
  
For any bad sequence $\langle q_0, q_1, \ldots, q_n\rangle$ in $Q$ we
know that $i<j\le n$ implies that either $q_i$ is incomparable with
$q_j$ or $q_i > q_j$ and hence, in the latter case
$\rho(q_i)>\rho(q_j)$. Fixing a bad sequence $s=\langle q_0, q_1,
\ldots, q_n\rangle$, consider the set
\begin{equation*}
S_{s}\eqdef\{ \langle q_{i_0}, q_{i_1}, \ldots, q_{i_m}\rangle :\, 
i_0<i_1\cdots<i_m=n
\land
 \rho(q_{i_0})\le  \rho(q_{i_1})\cdots\le
 \rho(q_{i_m})
\}.
\end{equation*}
In other words,
$S_s$ consists of subsequences of $s$ that end with $q_n$ and where
all elements are incomparable. So for each
$t\in S_{s}$ the value $g(t)$ is defined. We define that $\varphi(s)$
is the minimum over all $g(t)$ for $t\in S_{s}$. The intuition here is
that $\varphi$ is an ordinal measure for the longest incomparable
sequence within a bad sequence.  Now we are going to combine $\rho$
and $\varphi$ into a function $f$ defined on bad sequences. Given such
a sequence $s=\langle q_0, q_1, \ldots, q_n\rangle$, we let
\begin{equation*}
f(s)\eqdef \bigl\langle 
\bigl(\rho(q_0), \varphi(\langle q_0\rangle)\bigr),
\bigl(\rho(q_1), \varphi(\langle q_0, q_1\rangle)\bigr), 
\ldots,
\bigl(\rho(q_n), \varphi(\langle q_0, q_1,\ldots, q_n\rangle)\bigr)
\bigr\rangle.
\end{equation*}
Noticing that every non-empty subsequence of a bad sequence is bad, we
see that $f$ is a well-defined function which maps $\Bad(Q)$ into the
set of finite sequences from $\alpha\times\beta$. Moreover, let us
notice that every sequence in the image of $f$ is a bad sequence in
$\alpha\times\beta$: if $i<j$ 
and 
$\rho(q_i)\le \rho(q_j)$,  let $t$
be a sequence from $S_{\langle q_0, q_1, q_2, \ldots
q_i\rangle}$ such that $g(t)=\varphi(\langle q_0, q_1, q_2, \ldots
q_i\rangle)$. Hence $t$ includes $q_i$ and for every $q_k\in t$ we
have $\rho(q_k)\le \rho(q_i)\le \rho(q_j)$. Therefore $t\frown q_j$
was taken into account when calculating $\varphi(\langle q_0, q_1,
q_2, \ldots q_j\rangle)$. In particular,
\begin{equation}
\varphi(\langle q_0, q_1, q_2, \ldots q_j\rangle)\le g(t\frown q_j)< g(t)=\varphi(\langle q_0, q_1, q_2, \ldots q_i\rangle)\:.
\end{equation}
Then
 $(\rho(q_i), \varphi(\langle q_0, q_1, q_2, \ldots
q_i\rangle))\not\le (\rho(q_j), \varphi(\langle q_0, q_1, q_2, \ldots
q_j\rangle))$. 
Another possibility when $i<j$ is that $\rho(q_i)>\rho(q_j)$ and it
yields the same conclusion.
We have therefore shown that $f:\Bad(Q)\into
\Bad(\alpha\times\beta)$. Let us also convince ourselves that $f$ is a
tree homomorphism, meaning a function that preserves the strict tree
order.  The tree $\Bad(Q)$ is ordered by  initial segments, the
order which we have denoted by $\initial$. If $s\initial t$, then obviously $f(s)\initial f(t)$.  Given that it is well
known and easy to see that tree homomorphisms can only increase the
rank of a tree, we have that $\o(Q)\le \o(\alpha\times\beta)$.  The
latter, as shown by \citet{deJonghParikh}, is equal to
$\alpha\otimes\beta=\w(Q)\otimes \h(Q)$ (note that $\otimes$ is
commutative).
\ifams\relax\else\hfill\qedsymbol\fi
\end{proof}

From Thm.~\ref{thm-oandh} we derive a useful consequence.  Recall that
$\alpha$ is \emph{additive} (or \emph{multiplicative})
\emph{principal} if $\beta,\gamma<\alpha$ implies
$\beta+\gamma<\alpha$ (respectively implies $\beta \cdot
\gamma<\alpha$).  These implications also hold for natural sums and
products.\todo[size=\tiny]{I've put a proof in the Appendix for our
  peace of mind}
\begin{corollary}
\label{thm-w=o-4multprinc}
Assume that $\o(Q)$ is a principal multiplicative ordinal and that
$\h(Q) < \o(Q)$.  Then $\w(Q) = \o(Q)$.
\end{corollary}
\begin{proof}\renewcommand{\qedsymbol}{$\eop_{\ref{thm-w=o-4multprinc}}$}
Assume, by way of contradiction, that $\w(Q)<\o(Q)$.  From
$\h(Q)<\o(Q)$ we deduce $\h(Q)\otimes\w(Q)<\o(Q)$ (since $\o(Q)$ is
multiplicative principal), contradicting the inequality \eqref{eq-KT-ineq} in
Thm.~\ref{thm-oandh}. Hence $\w(Q) \geq \o(Q)$, and necessarily $\w(Q) =
\o(Q)$, again by~\eqref{eq-KT-ineq}.~%
\ifams\relax\else\hfill\qedsymbol\fi
\end{proof}

\section{Computing the Invariants of Common WQOs}
\label{sec-computing-w}

We now consider WQOs obtained in various well-known ways and address
the question of computing their width, and recall along the way what
is known about their height and maximal order type.

In the ideal case, there would be a means of defining
well-quasi-orders as the closure of some simple orders, in the
`Hausdorff-like' spirit of Thm.~\ref{allFAC}.  Unfortunately, no such
result is known and indeed it is unclear which class of orders one
could use as a base---for example how would one obtain Rado's example
(see Sec.~\ref{sec-Rado}) from a base of any `reasonable orders.'
Therefore, our study of the width of WQO orders will have to be
somewhat pedestrian, concentrating on concrete situations.
\ifams\relax\else\vspace{-.5em}\fi
\subsection{Lexicographic Sums}
In the case of lexicographic sums along an ordinal (defined in
Sec.~\ref{wqoversusall}), we have the following result.

\begin{lemma}\label{theorem-lexsum}
Suppose that for an ordinal $\alpha$ we have a family of WQOs
$\{P_i:\, i<\alpha\}$.  Then $\Sigma_{i<\alpha}P_i$ is a WQO, and:
\begin{enumerate}
\item 
$\o(\Sigma_{i<\alpha}P_i)=\Sigma_{i<\alpha} \o(P_i)$,
\item 
$\h(\Sigma_{i<\alpha}P_i)=\Sigma_{i<\alpha} \h(P_i)$,
\item 
$\w(\Sigma_{i<\alpha}P_i)=\sup_{i<\alpha} \w(P_i)$.
\end{enumerate}
\end{lemma}

\begin{proof} 
First note that any infinite bad sequence in $\Sigma_{i<\alpha}P_i$
would either have an infinite projection to $\alpha$ or an infinite
projection to some $P_i$, which is impossible.  Hence
$\Sigma_{i<\alpha}P_i$ is a WQO. Therefore the values
$\w(\Sigma_{i<\alpha}P_i)$, $\o(\Sigma_{i<\alpha}P_i)$ and
$\h(\Sigma_{i<\alpha}P_i)$ are well defined.
\begin{enumerate}
\item 
We use Thm.~\ref{thm-equivalences-2}. Let
$\alpha_i\eqdef\o(P_i)$, then $\Sigma_{i<\alpha}\alpha_i$ is isomorphic to
a linearisation of $\Sigma_{i<\alpha}P_i$. Hence
$\o(\Sigma_{i<\alpha}P_i) \ge \Sigma_{i<\alpha} \o(P_i)$. Suppose that
$L$ is a linearisation of $\Sigma_{i<\alpha}P_i$ (necessarily a well
order), then the projection of $L$ to each $P_i$ is a linearisation of
$P_i$ and hence it has type $\le \alpha_i$. This gives that the type
of $L$ is $\le \Sigma_{i<\alpha} \alpha_i$, proving the other side of
the desired inequality.

\item 
We use Thm.~\ref{thm-equivalences-1}. Any chain
$C$ in $\Sigma_{i<\alpha}P_i$ can be obtained as
$C=\Sigma_{i<\alpha}C_i$, where $C_i$ is the projection of $C$ on the
coordinate $i$.  The conclusion follows as in the case of $\o$.
\item
Every non-empty sequence of incomparable elements in $P$ must come
from one and only one $P_i$, hence $\Inc(P)=\bigcup_{i\in L}
\Inc(P_i)$, and therefore $\w(P_i)=\sup_{i<\alpha}\w(P_i)$ by Lem.~\ref{bunching}.
$\eop_{\ref{theorem-lexsum}}$
\end{enumerate}
\end{proof}

\subsection{Disjoint Sums}

We also defined disjoint sums in Sec.~\ref{wqoversusall} as sums along
an antichain.

\begin{lemma}\label{ABresults-disjsum} Suppose that $P_1,P_2,\ldots$
is a family of WQOs.
\begin{enumerate}
\item $\o(P_1\sqcup P_2) = \o(P_1)\oplus\o(P_2)$,
\item $\h(\bigsqcup_i P_i) = \sup \{\h(P_i)\}_i$,
\item $\w(P_1\sqcup P_2)=\w(P_1)\oplus \w(P_2)$.
\end{enumerate}
\end{lemma}\ifams\relax\else\clearpage\fi

\begin{proof}\renewcommand{\qedsymbol}{$\eop_{\ref{ABresults-disjsum}}$}
(1) is Thm.~3.4 from~\cite{deJonghParikh}.

\noindent
(2) is clear since, for an arbitrary family $P_i$ of WQOs,
$\Dec(\bigcup_i P_i)$ is isomorphic to $\bigsqcup_i \Dec(P_i)$. We
observe that, for infinite families, $\bigsqcup_i P_i$ is not WQO, but
it is still well-founded hence has a well-defined height.

\noindent
(3) is Lem.~1.10 from~\cite{AbBo} about antichain rank, which
translates to widths thanks to Thm.~\ref{equal}.%
\ifams\relax\else\qedsymbol\fi
\end{proof}

We can apply lexicographic sums to obtain the existence of WQO posets
of every width.

\begin{corollary}\label{theorem-obtained} For every ordinal $\alpha$,
  there is a WQO poset $P_\alpha$ such that $\w(P_\alpha)=\alpha$. 
\end{corollary}

\begin{proof} 
The proof is by induction on $\alpha$.  For $\alpha$ finite, the
conclusion is exemplified by an antichain of length $\alpha$. For
$\alpha$ a limit ordinal let us fix for each $\beta<\alpha$ a WPO
$P_\beta$ satisfying $\w(P_\beta)=\beta$.  Then
$\w(\Sigma_{\beta<\alpha} P_\beta)=\sup_{\beta<\alpha} \beta=\alpha$,
as follows by Lem.~\ref{theorem-lexsum}. For $\alpha=\beta+1$, we take
$P_\alpha=P_\beta\sqcup 1$, i.e., $P_\beta$ with an extra
(incomparable) element added, and rely on $\w(Q\sqcup 1)=\w(Q)\oplus
1=\w(Q)+1$ shown in Lem.~\ref{ABresults-disjsum}.
$\eop_{\ref{theorem-obtained}}$
\end{proof}

\subsection{Direct Products}
Direct products are again a particular case of lexicographic sums
along a poset~$Q$, this time of the same poset $P$.  While the cases
of $\o$ and $\h$ are mostly folklore, the width of $P\cdot Q$ is not
so easily understood, and its computation in Lem.~1.11 from
\cite{AbBo} uses the notion of \emph{Heisenberg products} $\alpha\odot
\beta$, defined for any ordinal $\alpha$ by induction on the ordinal
$\beta$:
\begin{align*}
  \alpha\odot 0&\eqdef 0\:,&
  \alpha\odot (\beta+1)&\eqdef (\alpha\odot \beta)\oplus \alpha\:,&
  \alpha \odot \lambda&\eqdef\sup\{(\alpha\odot
  \gamma)+1:\,\gamma<\lambda\}
\end{align*}
where $\lambda$ is a limit ordinal.  Note that this differs from the
natural product, and is not commutative: $2\odot \omega=\omega$ but
$\omega\odot 2=\omega\cdot 2$.

\begin{lemma}[\citeauthor{AbBo}]\label{ABresults-directprod} Suppose
  that $P$ and $Q$ are two WPOs.
\begin{enumerate}
\item $\o(P\cdot Q)=\o(P)\cdot \o(Q)$,
\item $\h(P\cdot Q)=\h(P)\cdot \h(Q)$,
\item $\w(P\cdot Q)=\w(P)\odot \w(Q)$.
\end{enumerate}
\end{lemma}

\subsection{Cartesian Products}\label{finiteproducts}

The next simplest operation on WQOs is their Cartesian product.  It
turns out that the simplicity of the operation is deceptive and that 
the height and, especially, the width of a product $P\times Q$ are not
as simple as we would like. As a consequence, this section only
provides partial results and is unexpectedly long.

To
recall, the product order $P\times Q$ of two partial orders is defined
on the pairs $(p,q)$ with $p\in P$ and $q\in Q$ so that $(p,q)\le
(p',q')$ iff $p\le _P p'$ and $q\le_Q q'$. It is easy to check and
well known that product of WQOs is WQO and similarly for FAC and WF
orders.

The formula for calculating $\o(P\times Q)$ is still simple. It was first established by
\citet[Thm.~3.5]{deJonghParikh}; see also \citep[Thm.~6]{BlGu}.
\begin{lemma}[\citeauthor{deJonghParikh}]\label{th-oprod}
  Suppose that $P$ and $Q$ are two WQOs.  Then $\o(P\times
  Q)=\o(P)\otimes\o(Q)$.
\end{lemma}

The question of the height of products is also well studied and a
complete answer appears in \citep{Abraham-Dilworth}, where it is
stated that the theorem is well known. The following statement is a
reformulation of Lem.~1.8 of \cite{Abraham-Dilworth}.

\begin{lemma}[Abraham; folklore]\label{the-heightproducts} 
If $\rho_P:\, P\to \h(P)$ and $\rho_Q:\, Q\to \h(Q)$ are the rank
functions of the well-founded posets $P$ and $Q$, then the rank
function $\rho$ on $P \times Q$ is given by $\rho (x, y) = \rho_P (x)
\oplus \rho_Q (y)$.  In particular,
\[ 
 \h(P \times  Q)=\sup \{\alpha\oplus\beta+1 :\, \alpha<\h(P)\land
 \beta<\h(Q)\}
\:.
\]
\end{lemma}

We recall that for any two ordinals $\alpha$ and $\beta$ we have
$\sup_{\alpha'<\alpha, \beta'<\beta} \alpha'\oplus \beta'+1 <
\alpha\oplus \beta$ \citep[see e.g.][p.~55]{AbBo}, thus the statement
in Thm.~\ref{the-heightproducts} cannot be easily simplified.

\begin{remark}[Height of products of finite ordinals]
\label{lem-h-nxm}
The very nice general proof of \citet[Lem.~1.8]{Abraham-Dilworth} can
be done in an even more visual way in the case of finite ordinals.  Let
$P=n_1\times \cdots \times n_k$ for some finite
$n_1,\ldots,n_k\in\omega$; then $\h(P)=n_1+\cdots+n_k+1-k$.

Indeed, we observe that any chain $\mathbf{a}_1 <_P\cdots
<_P\mathbf{a}_\ell$ in $P$ leads to a strictly increasing
$|\mathbf{a}_1|<\cdots<|\mathbf{a}_\ell|$, where by $|{\mathbf a}|$ we
denote the sum of the numbers in ${\mathbf a}$. Since
$|\mathbf{a}_\ell|$ is at most $\sum_i (n_i-1)=(\sum_i n_i)-k$ and
since $|\mathbf{a}_1|$ is at least $0$, the longest chain has length
$1+\sum_i n_i-k$. Furthermore it is easy to build a witness for this length. We conclude by
invoking Thm.~\ref{thm-equivalences-1} which states that for any WPO
$P$, $\h(P)$ is the length of the longest chain in $P$.
$\eop_{\ref{lem-h-nxm}}$
\end{remark}

Having dealt with $\h$ and $\o$, we are left with $\w$. Here we cannot
hope to have a uniform formula expressing $\w(P\times Q)$ as a
function of $\w(P)$ and $\w(Q)$.  Indeed, already in the case of
ordinals one always has $\w(\alpha)=\w(\beta)=1$, while
$\w(\alpha\times\beta)$ has quite a complex form, as we are going to
see next.

\subsubsection{Products of Ordinals}\label{prod-ordinals}
Probably the simplest example of WQO which is not actually an ordinal,
is provided by the product of two ordinals.  Thanks to
Thm.~\ref{equal}, we can translate results of \cite{Abraham-Dilworth},
Section 3 to give a recursive formula which completely characterises
$\w(\alpha\times \beta)$ for $\alpha,\beta$ ordinals. We shall sketch
how this is done.

First note that if one of $\alpha,\beta$ is a finite ordinals $n$, say $\alpha=n$, then we have $\w(n \times \beta)=\min\{n,\beta\}$. The next
case to consider is that of successor ordinals, which is taken care by the following Thm.~\ref{successorcase}. Abraham proved this theorem using
the method of residuals and induction, we offer an alternative proof using the rank of the tree $\Inc$.

\begin{theorem}[Abraham]\label{successorcase} For any ordinals
$\alpha,\beta$ with $\alpha$ infinite, we have $\w(\alpha\times (\beta+1))=\w(\alpha\times \beta)+1$.
\end{theorem}
The proof is provided by the next two lemmas.
\begin{lemma}
\label{lem2-axb+1}
$\w(\alpha\times (\beta+1))\leq \w(\alpha\times\beta)+1$ for any
ordinals $\alpha,\beta$.
\end{lemma}

\begin{proof}\renewcommand{\qedsymbol}{$\eop_{\ref{lem2-axb+1}}$}
Write $I$ for $\Inc(\alpha\times (\beta+1))$ and $I'$ for
$\Inc(\alpha\times \beta)$. Any sequence $s = \langle p_1,\ldots,
p_\ell\rangle$ which is in $I$, is either in $I'$ or contains a single pair of the
form $p_i = (a,\beta)$, with $a<\alpha$. In the latter case we write $s'$
for $s$ with $p_i$ removed. Note that $s'$ is in $I'$ (except when $s$
has length 1). Let $\rho':I'\to \rank(I')=\w(\alpha\times\beta)$ be a
ranking function for $I'$ and define $\rho:I\to \ON$ via
\[
  \rho(s) \egdef\begin{cases}
  \rho'(s)+1 & \text{if $s\in I'$,}\\
  \rho'(s') & \text{if $s\not\in I'$ and $|s|>1$,}\\
  \rank(I') & \text{otherwise.}
  \end{cases}
\]
One easily checks that $\rho$ is anti-monotone. For this assume $s
\initial t$: (1) if both $s$ and $t$ are in $I'$, monotonicity is
inherited from $\rho'$; (2) if none are in $I'$ then
$s'\initial t'$ (or $s'$ is empty) and again monotonicity is
inherited (or $\rho(s) = \rank(I') > \rho'(t') =
\rho(t)$); (3) if $s$ is in $I'$ and $t$ is not then $s\initialeq
t'$, entailing $\rho'(s) \geq \rho'(t')$ so that $\rho(s) =
\rho'(s)+1 > \rho'(t') = \rho(t)$.

In conclusion $\rho$, having values in $\w(\alpha\times\beta)+1$,
witnesses the assertion of the lemma.
\ifams\relax\else\qedsymbol\fi
\end{proof}

\begin{lemma}
\label{lem1-axb+1}
If $\alpha$ is infinite then $\w(\alpha\times (\beta+1))\geq
\w(\alpha\times\beta)+1$ for any $\beta$.
\end{lemma}
\begin{proof}\renewcommand{\qedsymbol}{$\eop_{\ref{lem1-axb+1}}$}
Write $I$ for $\Inc(\alpha\times\beta)$. Any $s\in I$
has the form $s=\langle(a_1,b_1),\ldots,(a_\ell,b_\ell)\rangle$. We
write $s_+$ for the sequence $\langle(a_1+1,b_1), \ldots,
(a_\ell+1,b_\ell)\rangle$ and observe that it is still a sequence over
$\alpha\times\beta$ since $\alpha$ is infinite, and that its elements
form an antichain (since the elements of $s$ did). Let now $s'_+$ be
$r\frown s_+$ where $r=\langle(0,\beta)\rangle$: the prepended
element is not comparable with any element of $s_+$ so that $s'_+$ is
an antichain and $s'_+\initialeq t'_+$ iff $s_+\initialeq t_+$ iff
$s\initialeq t$. Write $I'_+$ for $\{s'_+~|~s\in I\}\cup \{r\}$. This
is a tree made of a root glued below a tree isomorphic to $I$. Hence
$\rank(I'_+)=\rank(I)+1$. On the other hand, $I'_+$ is a substructure
of $\Inc(\alpha\times(\beta+1))$ hence
$\w(\alpha\times(\beta+1))\geq\rank(I'_+)$.
\ifams\relax\else\qedsymbol\fi
\end{proof}

With Thm.~\ref{successorcase} in hand, the remaining case is to compute $\w(\alpha\times\beta)$ when $\alpha, \beta$ are limit ordinals.
This translates into saying that $\alpha=\omega\alpha'$ and $\beta=\omega\beta'$ for some $\alpha', \beta'>0$. A recursive formula describing the
weight of this product is the main theorem 
of Section 3 of \cite{Abraham-Dilworth}, which we now quote. It is proved using a complex application of the method of residuals and induction.

\begin{theorem}[Abraham]\label{Abrahamlimits} Suppose that $\alpha$ and $\beta$ are given in their Cantor normal forms
$\alpha=\omega^{\alpha_0}\cdot m_0+\rho$, $\beta=\omega^{\beta_0}\cdot
n_0+\sigma$, where $\omega^{\alpha_0}\cdot m_0$ and
$\omega^{\beta_0}\cdot n_0$ are the leading terms and
$\rho$ and $\sigma$ are the remaining terms of the Cantor normal forms of $\alpha$
and $\beta$ respectively.
Then if $\alpha=1$, we have $\w(\omega\times\omega\beta)=\omega\beta$, and in general
\[
\w(\omega\alpha\times \omega\beta)=
\omega\omega^{\alpha_0\oplus\beta_0}\cdot(m_0+n_0-1) \oplus 
\w(\omega\omega^{\alpha_0}\times \omega\sigma)\oplus \w(\omega\omega^{\beta_0}\times \omega\rho).
\]
\end{theorem}

It would be interesting to have a closed rather than a recursive formula for the width of the product of two ordinals. However, the formula does give us a closed form of values of the weight
of the product of two ordinals with only one term in the Cantor normal form, as we now remark. Here $m,n$ are finite ordinals $\ge 1$.

\begin{enumerate}
\item If $k,\ell<\omega$ then we have 
\[\w(\omega^{1+k}\cdot m \times
  \omega^{1+\ell}\cdot n)=
\w\bigl(\omega(\omega^k\cdot m)\times\omega(\omega^\ell\cdot
n)\bigr)=\omega^{k+\ell-1}\cdot(m+n-1)\:.
\]
\item (example 3.4 (3) from \cite{Abraham-Dilworth}) If $\alpha,\beta\ge\omega$ then $1+\alpha=\alpha$ and $1+\beta=\beta$, so 
\[
\w(\omega^\alpha\cdot m\times \omega^\beta\cdot n)=
\w\bigl(\omega(\omega^\alpha\cdot m)\times \omega(\omega^\beta\cdot n)\bigr)=
\omega^{\alpha\oplus\beta}\cdot(m+n-1)\:.
\]
\item If $\alpha\ge\omega$ and $k<\omega$ then $\w(\omega^\alpha\cdot m\times \omega^{1+k}\cdot n)=\omega^{\alpha +k}\cdot(m+n-1)$.
\end{enumerate}

Let us mention one more result derivable from Thm.~\ref{Abrahamlimits}.
\begin{lemma}[Abraham]
\label{lem-w-omega-x-alpha}
$\w(\omega\times \alpha)=\alpha$  for any ordinal $\alpha$.
\end{lemma}
\begin{proof}\renewcommand{\qedsymbol}{$\eop_{\ref{lem-w-omega-x-alpha}}$}
By induction on $\alpha$. If $\alpha$ is a limit, we write it
$\alpha=\omega\alpha' = \omega(\omega^{\alpha_0}\cdot m_0+\cdots+
\omega^{\alpha_\ell}\cdot m_\ell)$. Now Thm.~\ref{Abrahamlimits}
yields $\w(\omega\times\omega\alpha')= \omega\omega^{\alpha_0}\cdot
m_0 \oplus \cdots\oplus \omega\omega^{\alpha_\ell}\cdot
m_\ell=\alpha$. If $\alpha$ is a successor, we use
Lem.~\ref{lem2-axb+1} and~\ref{lem1-axb+1}.
\ifams\relax\else\qedsymbol\fi
\end{proof}

\subsubsection{Finite Products and Transferable Orders}\label{gamesonproduct}

Since the width of the product of two ordinals is understood, we can
approach the general question of the width of products of two or a
finite number of WQO posets $P_i$ by reducing it to the width of some
product of ordinals. Using that strategy, we give a lower bound to
$\w(\Pi_{i\le n}P_i)$.

\begin{theorem}\label{thm-lboundPxQ} 
For any WQO posets $P_0, P_1\ldots P_n$, $\w(\prod_{i\le n}P_i )\ge
\w(\prod_{i\le n}\h(P_i))$.
\end{theorem}
The proof follows directly from a simple lemma, which is of
independent interest:

\begin{lemma}\label{lem-htintoproduct} 
Suppose that $P_0, P_1\ldots P_n$ are WQO posets. Then $\prod_{i\le
n}\h(P_i)$ embeds into $\prod_{i\le n}P_i $ as a substructure.
\end{lemma}

\begin{proof}\renewcommand{\qedsymbol}{$\eop_{\ref{lem-htintoproduct}}$}
We use Thm.~\ref{thm-equivalences-2} and pick, in each
$P_i$, a chain $C_i$ in $P_i$ that has
order type $\h(C_i)=\h(P_i)$.  Then $\prod_{i\le n}C_i$
is an induced suborder of $\prod_{i\le n}\h(P_i)$ which is isomorphic
to $\prod_{i\le n}\h(P_i)$.
\ifams\relax\else\qedsymbol\fi
\end{proof}

Now we shall isolate a special class of orders for which it will be
possible to calculate certain widths of products.  Let us write $\down
x$ for the downwards-closure of an element $x$, i.e., for $\{y: x\leq
y\}$.

\begin{definition}\label{def-everywheredense} A FAC partial order $P$ belongs to the class $\TT$ of \emph{transferable orders} if 
 $\w(P\setminus (\down x_1 \cup\cdots\cup\down x_n))=\w(P)$ for any
(finitely many) elements $x_1,\ldots,x_n\in P$.
\end{definition}

\begin{theorem}\label{thm-foisk} Suppose that $P$ is a WQO
  transferable poset and $\delta$ is an ordinal. Then $\w(P\times
  \delta)\ge \w(P)\cdot \delta$.
\end{theorem}

\begin{proof}\renewcommand{\qedsymbol}{$\eop_{\ref{thm-foisk}}$}
Write $\gamma$ for $\w(P)$: we prove that Player 2 has a winning
strategy, denoted $\sigma_{P'\times\delta,\alpha}$, for each game
$G_{P'\times\delta,\alpha}$ where $P'$ is some $P\setminus(\down
y_1\cup\cdots\cup\down y_n)$ and $\alpha\leq \gamma\cdot\delta$.

The proof is by induction on $\delta$.

If $\delta=0$ then $\alpha=0$ and Player 1 loses immediately.

If $\delta=\lambda$ is a limit, the strategy for Player 2 depends on
Player 1's first move. Say it is
$\alpha'<\alpha\leq\gamma\cdot\delta$. Then
$\alpha'<\gamma\cdot\delta$ means that $\alpha'<\gamma\cdot\delta'$
for some $\delta'<\delta$. Player 2 chooses one such $\delta'$ and now
applies $\sigma_{P'\times\delta',\alpha'+1}$ (which exists and is
winning by the induction hypothesis) for the whole game. Note that a strategy for a
substructure $P'\times\delta'$ of the original $P'\times\delta$ will
lead to moves that are legal in the original game. Also note that
$\alpha'+1$ is $\leq\gamma\cdot\delta'$.

If $\delta=\epsilon+1$ is a successor then Player 2 answers each move
$\alpha_1,\ldots,\alpha_m$ played by Player 1 by writing it in the
form $\alpha_i=\gamma\cdot\delta_i+\beta_i$ with $\beta_i<\gamma$.
Note that $\delta_i<\delta$. If $\delta_1=\cdots=\delta_m=\epsilon$, note that $\beta_1>\beta_2>\ldots \beta_m$.
Let Player 2 play $(x_m,\epsilon)$ where $x_m$ is $\sigma_{P',\gamma}$
applied on $\beta_1,\ldots,\beta_m$ (that strategy exists and is
winning since $P$ is transferable and has width $\gamma$). If
$\delta_m<\epsilon$ then Player 2 switches strategy and now uses
$\sigma_{P''\times\epsilon,\gamma\cdot\epsilon}$ as if a new game was
starting with $\alpha_m$ as Player 1's first more, and for
$P''=P'\setminus(\down x_1\cup\cdots\cup\down x_{m-1})$. By
the induction hypothesis , Player 2 will win by producing a sequence $S''$ in
$P''\times\epsilon$. These moves are legal since
$(x_1,\epsilon)\cdots(x_{m-1},\epsilon)\frown S''$ is an antichain in
$P'\times(\epsilon+1)$.
\ifams\relax\else\qedsymbol\fi
\end{proof}

In order to use Thm.~\ref{thm-foisk}, we need actual instances of
transferable orders.
\begin{lemma}\label{lem-omegatransferable} For any $1\leq\alpha_1,\ldots,\alpha_n$, the order $P=\omega^{\alpha_1}\times\cdots\times \omega^{\alpha_n}$ is transferable.
\end{lemma}

\begin{proof}\renewcommand{\qedsymbol}{$\eop_{\ref{lem-omegatransferable}}$}
Since each $\omega^{\alpha_i}$ is additive principal, $P\setminus (\down
x_1\cup\cdots\cup\down x_m)$ contains an isomorphic copy of $P$ for
any finite sequence $x_1,\ldots,x_m$ of elements of $P$.
\ifams\relax\else\qedsymbol\fi
\end{proof}

\begin{theorem}\label{cor-productofsquares} Let $P$ be a transferable
WPO poset.
\begin{enumerate}
\item  Suppose that $1\le m<\omega$. Then  $\w(P) \cdot m\le \w(P\times m)\le \w(P) \otimes m$.

\item If $\w(P)=\omega^\gamma$ for some $\gamma$, then $\w(P\times m)=\w(P) \cdot m$ (Note that this applies to any $P$ which is the product of the
form $\omega^\alpha\times \omega^\beta$, see the examples after Thm.~\ref{Abrahamlimits}).

\item $\w(\omega\times\omega\times\omega)=\omega^2$.

\end{enumerate}
\end{theorem}

An easy way to provide an upper bound needed in the proof of Thm.~\ref{cor-productofsquares} is given by the following observation:

\begin{lemma}\label{lem-uboundPxQ} 
For any FAC poset $P$ and $1\le m<\omega$, $\w(P\times m)\le
\w(P)\otimes m$.
\end{lemma}

\begin{proof}\renewcommand{\qedsymbol}{$\eop_{\ref{lem-uboundPxQ}}$}
We just need to remark that $P\times m$ is an augmentation of the
perpendicular sum $\sqcup_{i < m}P$ and then apply Lem.~\ref{ABresults-disjsum}.
\ifams\relax\else\qedsymbol\fi
\end{proof}

\begin{proof}[\ifams Proof \fi of Thm.~\ref{cor-productofsquares}]
\renewcommand{\qedsymbol}{$\eop_{\ref{cor-productofsquares}}$}
(1) We get $\w(P\times m)\ge \w(P) \cdot m$ from Thm.~\ref{thm-foisk}.  
We get $\w(P\times m)\le \w(P) \otimes m$ from Lem.~\ref{lem-uboundPxQ}.

\noindent 
(2) This follows because $\omega^\gamma\otimes m = \omega^\gamma \cdot
m$.

\noindent 
(3) Let $P=\omega\times\omega$, hence we know that $\w(P)=\omega$.
Since any $P\times m$ is a substructure of $P\times \omega$, we
clearly have that $\w(P\times \omega)\ge \sup_{m<\omega} \w(P\times
m)= \sup_{m<\omega} \omega \cdot m=\omega^2$. Let us now give a proof
using games that $\w(P\times \omega)\le \omega^2$. It suffices to give
a winning strategy to Player 1 in the game $G_{P\times \omega,\gamma}$
for any ordinal $\gamma > \omega^2$.

So, given such a $\gamma$, Player 1 starts the game by choosing as his
first move the ordinal $\omega^2$. Player 2 has to answer by choosing
an element $x$ in $P\times \omega$, say an element $(p,m)$ with $p=(k,
\ell)$. Now notice that any element of $P\times \omega$ that is
incompatible with $(p,m)$ is either an element of $P\times m$ or of
the form $(q,n)$ for some $q\le p$ in $\omega\times \omega$, or is of
the form $(r,i)$ for some $r$ which is incompatible with $p$ in
$\omega\times \omega$. Therefore, any next step of Player 2 has to be
in an order $P'$ which is isomorphic to an augmentation of a
substructure of the disjoint union of the form
\begin{equation}\label{sqcups}
 P\times m\sqcup [(k+1)\times (\ell+1)] \times \omega \sqcup [(k+1)\times \omega]\times \omega \sqcup [(\ell +1)\times \omega]\times \omega.
\end{equation}
It now suffices for Player 1 to find an ordinal $o < \omega^2$
satisfying $o > \w(P')$ as the game will then be transferred to
$G_{P',o}$, where Player 1 has a winning strategy.  As $\omega^2$ is
closed under $\oplus$, it suffices to show that each of the orders
appearing in equation (\ref{sqcups}) has weight $<\omega^2$.  This is
the case for $P\times m$ by (2). We have that $\w\bigl( [(k+1)\times (\ell+1)]
\times \omega\bigr)= \w\bigl( (k+1)\times [(\ell+1) \times \omega]\bigr)$, which by
applying Lem.~\ref{lem-uboundPxQ} is $\le (\ell+1)\cdot (k+1)$. For
$[(k+1)\times \omega]\times \omega$, we apply 
Lem.~\ref{lem-uboundPxQ} to $\omega\times\omega$, to obtain
$\w\bigl([(k+1)\times \omega]\times \omega\bigr)\le \omega \cdot (k+1)$ and
similarly $\w\bigl([(\ell+1)\times \omega]\times \omega\bigr)\le \omega \cdot
(\ell+1)$.
\ifams\relax\else\qedsymbol\fi
\end{proof}

\subsection{Finite Multisets, Sequences, and Trees}\label{sec-sequences}

Well-quasi-orders are also preserved by building multisets, sequences,
and trees with WQO labels, together with suitable embedding relations.

\emph{Finite sequences} in $Q^{<\omega}$ are compared by the
\emph{subsequence embedding} ordering defined by
$s=\tup{x_0,\dots,x_{n-1}}\leq_* s'=\tup{x'_0,\dots,x'_{p-1}}$ if
there exists $f{:}\,n\to p$ strictly monotone such that $x_i \leq x'_{f(i)}$ in
$Q$ for all $i\in n$.  The fact that $(Q^{<\omega},{\leq_*})$ is WQO
when $Q$ is WQO was first shown by \citet{Higman}.

Given a WQO $(Q,{\leq})$, a \emph{finite multiset} over~$Q$ is a
function $m$ from $Q\to\+N$ with finite support, i.e.\ $m(x)>0$ for
finitely many $x\in Q$.  Equivalently, a finite multiset is a finite
sequence $m$ in $Q^{<\omega}$ where the order is irrelevant, and can
be noted as a `set with repetitions' $m=\{x_1,\dots,x_n\}$; we denote
by $M(Q)$ the set of finite multisets over~$Q$.  The \emph{multiset
  embedding} ordering is then defined by
$m=\{x_0,\dots,x_{n-1}\}\leqH m'=\{x'_0,\dots,x'_{p-1}\}$ if there
exists an injective function $f{:}\,n\to p$ with $x_i \leq x'_{f(i)}$
in $Q$ for all $i\in n$.  As a consequence of $(Q^{<\omega},{\leq_*})$
being WQO, $(M(Q),{\leqH})$ is also WQO when $Q$~is.

Finally, a (rooted, ordered) \emph{finite tree} $t$ over $Q$ is either
a leaf $x()$ for some $x\in Q$, or a term $x(t_1,\dots, t_n)$ for some
$n>0$, $x\in Q$, and $t_1,\dots,t_n$ trees over~$Q$.  
A tree has arity~$b$ if we bound $n$ by~$b$ in this definition.
We let $T(Q)$ denote the set of finite trees over~$Q$.  The
\emph{homeomorphic tree embedding} ordering is defined by
$t=x(t_1,\dots,t_n)\leqT t'=x'(t'_1,\dots,t'_p)$ (where $n,p\geq 0$)
if at least one the following cases occurs:
\begin{itemize}
\item $t\leqT t'_j$ for some $1\leq j\leq p$, or
\item $x\leq x'$ in $Q$ and $t_1\cdots t_n\leq_* t'_1\cdots t'_p$ for
  the subsequence embedding relation on $T(Q)$.
\end{itemize}
The fact that $(T(Q),{\leqT})$ is WQO when $Q$ is WQO was first shown
by \citet{Higman} for trees of bounded arity, before \citet{kruskal60}
proved it in the general case.  Note that it implies
$(Q^{<\omega},{\leq_*})$ being WQO for the special case of trees of
arity~$1$.

\subsubsection{Maximal Order Types}

The maximal order types of $M(Q)$, $Q^{<\omega}$, and $T(Q)$ have been
studied by \citet{weiermann2009} and \citet{schmidt79}; see also
\citet[Sec.~1.2]{vandermeeren15} for a nice exposition of these results.

For finite multisets with embedding, we need some additional
notations.  For an ordinal $\alpha$ with Cantor normal form
$\omega^{\alpha_1}+\cdots+\omega^{\alpha_n}$ where
$\o(P)\geq \alpha_1\geq \ldots\geq \alpha_n$, we let
\begin{equation}
  \widehat\alpha\eqdef\omega^{{\alpha_1}'}+\cdots+\omega^{{\alpha_n}'}
\end{equation}
where $\alpha'$ is $\alpha+1$ when $\alpha$ is an epsilon number,
i.e.\ when $\omega^\alpha=\alpha$, and is just $\alpha$ otherwise.

The following is \cite[Thm.~2]{weiermann2009}, with a corrected proof
due to \citet[Thm.~5]{VdMRaWe}.
\begin{theorem}[\citeauthor{weiermann2009}]\label{th-oM}
  Let $Q$ be a WQO.  Then $\o(M(Q))=\omega^{\widehat{\o(Q)}}$.
\end{theorem}
Thus, for $\o(Q)<\varepsilon_0$, one has simply
$\o(M(Q))=\omega^{\o(Q)}$.

\medskip
For finite sequences with subsequence embedding, we recall the following
result by \cite{schmidt79}.
\begin{theorem}[\citeauthor{schmidt79}]\label{th-oS}
  Let $Q$ be a WQO.  Then
\begin{equation*}
\label{eq-o-of-seq}
\o(Q^{<\omega})
=
\begin{cases}
  \omega^{\omega^{\o(Q)-1}} &\text{if $\o(Q)$ is finite},
  \\
  \omega^{\omega^{\o(Q)+1}} &\text{if $\o(Q)=\varepsilon+n$ for
    $\varepsilon$ an epsilon number and $n$ finite},
  \\
  \omega^{\omega^{\o(Q)}} &\text{otherwise}.
\end{cases}
\end{equation*}
\end{theorem}

The case of finite trees is actually a particular case of the results
of \citet{schmidt79} on embeddings in structured trees.  Her results
were originally stated using Sch\"utte's Klammer symbols, but can be
translated in terms of the $\vartheta$ functions of \citet{RaWe}.
Defining such ordinal notation systems is beyond the scope of this
chapter; it suffices to say for our results that the ordinals at hand
are going to be principal multiplicative.
\begin{theorem}[\citeauthor{schmidt79}]\label{th-oT}
  Let $Q$ be a WQO.  Then
  $o(T(Q))=\vartheta(\Omega^\omega\cdot\o(Q))$.
\end{theorem}

\subsubsection{Heights}\label{ssec-hast}

For a WQO $Q$ we define $\h^*(Q)$ as
\begin{equation}
\h^*(Q)\eqdef\begin{cases}
        \h(Q) & \text{if $\h(Q)$ is additive principal $\geq \omega$,}\\
        \h(Q)\cdot \omega & \text{otherwise.}
\end{cases}
\end{equation}

We are going to show that the heights of finite multisets, finite
sequences, and finite trees over $Q$ is the same, namely $\h^\ast(Q)$.

\begin{theorem}\label{eq-hT=hM=h*}
  Let $Q$ be a WF poset. Then
  $\h(M(Q))=\h(Q^{<\omega})=\h(T(Q))=\h^*(Q)$.
\end{theorem}
Since obviously $\h(M(Q))\leq \h(Q^{<\omega})\leq\h(T(Q))$, the claim
is a consequence of lemmata~\ref{lem-bound-hTQ}
and~\ref{lem-bound-hMQ} below.

\begin{lemma}\label{lem-bound-hTQ}
$\h(T(Q)) \leq \h^*(Q)$.
\end{lemma}
\begin{proof}\renewcommand{\qedsymbol}{$\eop_{\ref{cor-productofsquares}}$}
Consider a strictly decreasing sequence $x_0 >_T x_1 >_T \ldots$ in
$T(Q)$, where each $x_i$ is a finite tree over $Q$.  Necessarily these
finite trees have a nonincreasing number of nodes: $|x_0|\geq
|x_1|\geq\ldots$.  If we add a new minimal element $\bot$ below $Q$,
we can transform any $x_i$ by padding it with some $\bot$'s so that
now the resulting $x'_i$ has the same shape and size as $x_0$. Let us
use $1+Q$ instead of $\{\bot\}+Q$ so that the new trees belong to
$T(1+Q)$, have all the same shape, and form a strictly decreasing
sequence. This construction is in fact an order-reflection from
$\Dec(T(Q))$ to $\Dec\bigl(\bigsqcup_{n<\omega}(1+Q)^n\bigr)$, from
which we get
\begin{equation}
\label{eq-bound-TQ1}
\h(T(Q))\leq
\h(\bigsqcup_{n<\omega}(1+Q)^n)=\sup_{n<\omega} \h([1+Q]^n)
\:,
\end{equation}
using Lem.~\ref{ABresults-disjsum}.(2) for the last equality.  For
$n<\omega$, one has
\begin{equation}
\label{eq-bound-TQ2}
\h([1+Q]^n)
=\sup \{ (\alpha\otimes n)+1 ~:~ \alpha < 1+\h(Q) \}\:,
\end{equation}
using lemmata~\ref{theorem-lexsum}.(2)
and~\ref{the-heightproducts}.

If $\h(Q)\leq 1$, $\h(T(Q))=\h(Q)\cdot\omega=\h^*(Q)$ obviously.

For $\h(Q)> 1$, and thanks to \eqref{eq-bound-TQ1} and
\eqref{eq-bound-TQ2}, it is sufficient to show that $\alpha\otimes
n+1\leq \h^*(Q)$ for all $n<\omega$ and all $\alpha<1+\h(Q)$. We
consider two cases:
\begin{enumerate}
\item
If $\h(Q)\geq\omega$ is additive principal, $\alpha<1+\h(Q)=\h(Q)$
entails $\alpha\otimes n<\h(Q)$ thus $\alpha\otimes
n+1<\h(Q)=\h^*(Q)$.
\item
Otherwise the CNF for $\h(Q)$ is $\sum_{i=1}^m\omega^{\alpha_i}$ with
$m>1$.  Then $\alpha<1+\h(Q)$ implies $\alpha\leq
\omega^{\alpha_1}\cdot m$, thus $\alpha\otimes n+1\leq
\omega^{\alpha_1}\cdot m\cdot n +1\leq \omega^{\alpha_1+1}=\h(Q)\cdot
\omega=\h^*(Q)$.
\ifams\qedhere\else\qedsymbol\fi
\end{enumerate}
\end{proof}

\newcommand{\vx}{{\bm{x}}} \newcommand{\vy}{{\bm{y}}} Let us write
$M_n(Q)$ for the restriction of $M(Q)$ to multisets of size $n$.
\begin{lemma}\label{lem-MnQ-vs-Q^n}
$\h(M_n(Q))\geq \h(Q^n)$.
\end{lemma}
\begin{proof}\renewcommand{\qedsymbol}{$\eop_{\ref{lem-bound-hTQ}}$}
With $\vx=\tup{x_1,\ldots,x_n}\in Q^n$ we associate the multiset
$M_\vx = \{x_1,\ldots,x_n\}$. Obviously $\vx <_\times\vy$ implies
$M_\vx\leqH M_\vy$. We further claim that $M_\vy\not\leqH M_\vx$.
Indeed, assume by way of contradiction that $M_\vy\leqH M_\vx$. Then
there is a permutation $f$ of $\{1,\ldots,n\}$ such that $y_i\leq_Q
x_{f(i)}$ for all $i=1,\ldots,n$.  From $\vx\leq_\times\vy$, we get
\[
x_i\leq_Q y_i\leq_Q x_{f(i)}\leq_Q y_{f(i)}\leq x_{f(f(i))}
\leq y_{f(f(i))} \leq_Q \cdots \leq_Q x_{f^k(i)}\leq_Q y_{f^k(i)} \leq_Q\cdots
\]
So that for all $j$ in the $f$-orbit of $i$, $x_j\equiv_Q x_i\equiv_Q
y_j$, entailing $\vy\equiv_\times \vx$ which contradicts the
assumption $\vx<_\times \vy$.

We have thus exhibited a mapping from $Q^n$ to
$M_n(Q)$ that will map chains to chains. Hence $\h(Q^n)\leq
\h(M_n(Q))$.
\ifams\relax\else\qedsymbol\fi
\end{proof}

\begin{lemma}\label{lem-bound-hMQ}
$\h(M(Q)) \geq \h^*(Q)$.
\end{lemma}
\begin{proof}
The result is clear in cases where $\h^*(Q)=\h(Q)$ and when $\h(Q)=1$
entailing $\h(M(Q))=\omega=\h^*(Q)$. So let us assume that $\h(Q)$ is
not additive principal and has a CNF $\sum_{i=1}^m\omega^{\alpha_i}$
with $m>1$. Thus $\h^*(Q) = \h(Q)\cdot\omega =
\omega^{\alpha_1+1}$. Since by Lem.~\ref{the-heightproducts}, for
$0<n<\omega$, $\h(Q^n)=\sup \{\alpha\otimes n+1~:~\alpha<\h(Q)\}$, we
deduce $\h(Q^n)\geq \omega^{\alpha_1}\cdot n+1$%
.  Since $M_n(Q)$ is a
substructure of $M(Q)$, and using Lem.~\ref{lem-MnQ-vs-Q^n}, we
deduce
\begin{align*}
\hspace{1.9cm}\h(M(Q))&\geq \h(M_n(Q))\geq \h(Q^n)\geq \omega^{\alpha_1}\cdot n+1
\\
\shortintertext{for all $0<n<\omega$, hence}
\hspace{1.9cm}\h(M(Q))&\geq \sup_{n<\omega} \omega^{\alpha_1}\cdot n+1 =
\omega^{\alpha_1}\cdot \omega = \h^*(Q)\:.&\hspace{1.9cm}\eop_{\ref{lem-bound-hMQ}}
\end{align*}
\end{proof}

\subsubsection{Widths}
The previous analyses of the maximal order types and heights of
$M(Q)$, $Q^{<\omega}$, and $T(Q)$ allow us to apply the correspondence
between $\o$, $\h$, and $\w$ shown by \citet[Thm.~4.13]{kriz90b}, in
particular its consequence spelled out in Cor.~\ref{thm-w=o-4multprinc}.

\begin{theorem}
  \label{prop-w-M-seq-T}
  Let $Q$ be a WQO.  Then $\w(Q^\dagger)=\o(Q^\dagger)$ where
  $Q^\dagger$ can be $T(Q)$, or $Q^{<\omega}$ when $\o(Q)>1$, or
  $M(Q)$ when $\o(Q)>1$ is a principal additive ordinal.
\end{theorem}
\begin{proof}\renewcommand{\qedsymbol}{$\eop_{\ref{prop-w-M-seq-T}}$}
  First observe that $\h^\ast(Q)\leq
  \h(Q)\cdot\omega\leq\o(Q)\cdot\omega < \o(Q^\dagger)$ when
  $Q^\dagger$ is $T(Q)$ (by Thm.~\ref{th-oT}), $Q^{<\omega}$ with
  $\o(Q)>1$ (by Thm.~\ref{th-oS}), or $M(Q)$ with $\o(Q)>1$ (by
  Thm.~\ref{th-oM}).  Furthermore, when $Q^\dagger$ is $T(Q)$ or
  $Q^{<\omega}$, and when it is $M(Q)$ with $\o(Q)$ a principal
  additive ordinal, $\o(Q^\dagger)$ is a principal multiplicative
  ordinal.  Thus Cor.~\ref{thm-w=o-4multprinc} shows that
  $\w(Q^\dagger)=\o(Q^\dagger)$.
\ifams\relax\else\qedsymbol\fi
\end{proof}
The assumptions in Thm.~\ref{prop-w-M-seq-T} seem necessary.  For
instance, if $Q=1$, then $M(1)$ is isomorphic to $1^{<\omega}$ and
$\omega$, with height~$\omega$ and width~$1$.  If $A_3=1\sqcup 1\sqcup
1$ is an antichain with three elements, then $M(A_3)$ is isomorphic
with $\omega\times\omega\times\omega$, $\h(M(A_3))=\omega$ by
Lem.~\ref{the-heightproducts} or Thm.~\ref{eq-hT=hM=h*},
$\o(M(A_3))=\omega^3$ by Lem.~\ref{th-oprod}, and
$\w(M(A_3))=\omega^2$ by Thm.~\ref{cor-productofsquares}.(3).

\subsection{Infinite Products and Rado's Structure}\label{sec-Rado}

One may wonder what happens in the case of infinite products. We
remind the reader that the property of being WQO is in general not
preserved by infinite products. The classical example for this was
provided by Rado in \cite{rado54}, who defined what we call the
\emph{Rado structure}, denoted $(R,\leq)$: \footnote{We adopted the
  definition from~\cite{laver-wqos}.} Rado's order is given
as a structure on $\omega\times\omega$ where we define
\[
(a,b) \leq (a',b')\mbox{ if }[a=a'\mbox{ and }b\leq b']\mbox{ or }b<a'.
\]
The definition of BQOs was motivated by trying to find a property
stronger than WQO which is preserved by infinite products, so in
particular Rado's example is not a BQO \citep[see][Thm.~1.11 and
  2.22]{Milner}.

We can use the method of residuals and other tools
described in previous sections to compute.
\begin{xalignat}{3}
  \o(R) &= \omega^2, & \h(R)&=\omega, & \w(R)&=\omega,
\end{xalignat}
which gives the same ordinal invariants as those of the product $\omega\times\omega$, even though they are not isomorphic, and moreover
$\omega\times\omega$ is a BQO (since the notion of BQO is preserved under products) while Rado's order is not. Therefore one cannot characterise BQOs by the ordinal invariants considered here. Moreover, the two orders do not even embed into each other. To see this,
assume by way of contradiction that $f$ injects $\omega\times\omega$
into $R$. Write $(a_i,b_i)$ and $(c_i,d_i)$ for $f(0,i)$ and, resp.,
$f(i,0)$ when $i\in\omega$. Necessarily the $b_i$'s and the $d_i$'s
are unbounded. If the $a_i$'s are unbounded, one has the contradictory
$f(1,0)<_Rf(0,i)=(a_i,b_i)$ for some $i$, and there is a similar
contradiction if the $c_i$'s are unbounded, so assume the $a_i$'s and
the $c_i$'s are bounded by some $k$. By the pigeonhole principle, we
can find a pair $0<i,j$ with $a_i=c_j$ so that $f(0,i)\mathbin{\not\!\!\bot_R}
f(j,0)$, another contradiction. Hence $(\omega\times\omega)\not\leq R$.
In the other direction $R\not\leq(\omega\times\omega)$,
is obvious since $\omega\times\omega$ is BQO while $R$ is not.

\section{Concluding Remarks}\label{sec-concl}

We provide in Table~\ref{table-summary} a summary of our findings
regarding ordinal invariants of WQOs.  Mostly, the new results concern
the width $\w(P)$ of WQOs.  We note that the width $\w(P\times Q)$ of
Cartesian products is far from elucidated, the first difficulty being
that---unlike other constructs---it cannot be expressed as a function
of the widths $\w(P)$ and $\w(Q)$.  For Cartesian products,
Sec.~\ref{finiteproducts} only provide definite values for a few
special cases: for the rest, one can only provide upper and lower
bounds for the moment.

\begin{table}[tbp]
 \caption{Ordinal invariants of the main WQOs.
 \label{table-summary}}
 {\renewcommand{\arraystretch}{1.2}\ifams\small\setlength{\tabcolsep}{4pt}\else\setlength{\tabcolsep}{7pt}\fi
  \centering
  \begin{tabular}{cccc}
\toprule
$P$ & $\o(P)$ & $\h(P)$ & $\w(P)$
\\
\midrule
$\alpha\in\ON$ & $\alpha$ & $\alpha$ & $1$ (or $0$)
\\
$A_n$ (size $n$ antichain) & $n$ & $1$ & $n$
\\
Rado's $R$ & $\omega^2$ & $\omega$ & $\omega$
\\
\midrule
$\sum_{i\in\alpha}P_i$
  & $\sum_{i\in\alpha}\o(P_i)$ & $\sum_{i\in\alpha}\h(P_i)$ &
  $\sup_{i\in\alpha}\w(P_i)$
\\
$P\sqcup Q$ & $\o(P)\oplus \o(Q)$ & $\max(\h(P),\h(Q))$ & $\w(P)\oplus
\w(Q)$
\\
$P\cdot Q$ & $\o(P)\cdot \o(Q)$  & $\h(P)\cdot \h(Q)$
 & $\w(P)\odot \w(Q)$
\\
$P\times Q$ & $\o(P)\otimes \o(Q)$ &
$\sup_{\alpha<\h(P)\atop \beta< \h(Q)} \alpha\oplus\beta+1$ & see Sec.~\ref{finiteproducts}%
\\
\midrule
$M(P)$
&$\omega^{\widehat{\o(P)}}$&$\h^\ast(P)$, see Sec.~\ref{ssec-hast} & see Thm.~\ref{prop-w-M-seq-T}
\\
$P^{<\omega}$
&$\omega^{\omega^{\o(P)\pm 1}}$, see Thm.~\ref{th-oS}
&$\h^\ast(P)$ & $\o(P^{<\omega})$
\\
$T(P)$
&$\vartheta(\Omega^\omega\cdot\o(P))$ &$\h^\ast(P)$ & $\o(T(P))$
\\
\bottomrule
\end{tabular}}
\end{table}

\bibliographystyle{abbrvnat}
\bibliography{biblio}

\begin{thebibliography}{28}
\providecommand{\natexlab}[1]{#1}
\providecommand{\url}[1]{\texttt{#1}}
\expandafter\ifx\csname urlstyle\endcsname\relax
  \providecommand{\doi}[1]{doi: #1}\else
  \providecommand{\doi}{doi: \begingroup \urlstyle{rm}\Url}\fi

\bibitem[Abraham(1987)]{Abraham-Dilworth}
U.~Abraham.
\newblock A note on {Dilworth}'s theorem in the infinite case.
\newblock \emph{Order}, 4\penalty0 (2):\penalty0 107--125, 1987.
\newblock \doi{10.1007/BF00337691}.

\bibitem[Abraham and Bonnet(1999)]{AbBo}
U.~Abraham and R.~Bonnet.
\newblock Hausdorff's theorem for posets that satisfy the finite antichain
  property.
\newblock \emph{Fund. Math.}, 159\penalty0 (1):\penalty0 51--69, 1999.
\newblock ISSN 0016-2736.

\bibitem[Abraham et~al.(2012)Abraham, Bonnet, Cummings, D{\v{z}}amonja, and
  Thompson]{abcdzt}
U.~Abraham, R.~Bonnet, J.~Cummings, M.~D{\v{z}}amonja, and K.~Thompson.
\newblock A scattering of orders.
\newblock \emph{Trans. Amer. Math. Soc.}, 364\penalty0 (12):\penalty0
  6259--6278, 2012.
\newblock ISSN 0002-9947.
\newblock \doi{10.1090/S0002-9947-2012-05466-3}.

\bibitem[Blass and Gurevich(2008)]{BlGu}
A.~Blass and Y.~Gurevich.
\newblock Program termination and well partial orderings.
\newblock \emph{ACM Trans. Comput. Log.}, 9\penalty0 (3:18), 2008.
\newblock ISSN 1529-3785.
\newblock \doi{10.1145/1352582.1352586}.

\bibitem[de~Jongh and Parikh(1977)]{deJonghParikh}
D.~H.~J. de~Jongh and R.~Parikh.
\newblock Well-partial orderings and hierarchies.
\newblock \emph{Nederl. Akad. Wetensch. Proc. Ser. A {\bf 80}=Indag. Math.},
  39\penalty0 (3):\penalty0 195--207, 1977.
\newblock \doi{10.1016/1385-7258(77)90067-1}.

\bibitem[Dilworth(1950)]{Dilworth}
R.~P. Dilworth.
\newblock A decomposition theorem for partially ordered sets.
\newblock \emph{Ann. Math.}, 51\penalty0 (1):\penalty0 161--166, 1950.
\newblock \doi{10.2307/1969503}.

\bibitem[Dushnik and Miller(1941)]{DushnikMiller}
B.~Dushnik and E.~W. Miller.
\newblock Partially ordered sets.
\newblock \emph{Amer. J.~Math.}, 63\penalty0 (3):\penalty0 600--610, 1941.
\newblock ISSN 0002-9327.
\newblock \doi{10.2307/2371374}.

\bibitem[Higman(1952)]{Higman}
G.~Higman.
\newblock Ordering by divisibility in abstract algebras.
\newblock \emph{Proc. London Math. Soc. (3)}, 2:\penalty0 326--336, 1952.
\newblock ISSN 0024-6115.
\newblock \doi{10.1112/plms/s3-2.1.326}.

\bibitem[K{\v{r}}{\'{\i}}{\v{z}} and Thomas(1990)]{kriz90b}
I.~K{\v{r}}{\'{\i}}{\v{z}} and R.~Thomas.
\newblock Ordinal types in {Ramsey} theory and well-partial-ordering theory.
\newblock In J.~Ne{\v{s}}et{\v{r}}il and V.~R{\"{o}}dl, editors,
  \emph{Algorithms and Combinatorics vol.\ 5: Mathematics of Ramsey Theory},
  pages 57--95. Springer, 1990.
\newblock \doi{10.1007/978-3-642-72905-8_7}.

\bibitem[Kruskal(1960)]{kruskal60}
J.~Kruskal.
\newblock Well-quasi-ordering, the {T}ree {T}heorem, and {V}azsonyi's
  {C}onjecture.
\newblock \emph{Trans.~AMS}, 95\penalty0 (2):\penalty0 210--225, 1960.
\newblock \doi{10.2307/1993287}.

\bibitem[Kruskal(1972)]{Kruskal}
J.~B. Kruskal.
\newblock The theory of well-quasi-ordering: {A} frequently discovered concept.
\newblock \emph{J.~Comb. Theory A}, 13\penalty0 (3):\penalty0 297--305, 1972.
\newblock ISSN 0097-3165.
\newblock \doi{10.1016/0097-3165(72)90063-5}.

\bibitem[Laver(1976)]{laver-wqos}
R.~Laver.
\newblock Well-quasi-orderings and sets of finite sequences.
\newblock \emph{Math. Proc. Cambridge}, 79\penalty0 (1):\penalty0 1--10, 1976.
\newblock \doi{10.1017/S030500410005204X}.

\bibitem[Milner(1985)]{Milner}
E.~C. Milner.
\newblock Basic wqo- and bqo-theory.
\newblock In \emph{Graphs and order ({B}anff, {A}lta., 1984)}, volume 147 of
  \emph{NATO Adv. Sci. Inst. Ser. C Math. Phys. Sci.}, pages 487--502. Reidel,
  Dordrecht, 1985.
\newblock \doi{10.1007/978-94-009-5315-4_14}.

\bibitem[Milner and Sauer(1981)]{milner81}
E.~C. Milner and N.~Sauer.
\newblock On chains and antichains in well founded partially ordered sets.
\newblock \emph{J.~London Math. Soc.}, s2-24\penalty0 (1):\penalty0 15--33,
  1981.
\newblock \doi{10.1112/jlms/s2-24.1.15}.

\bibitem[Mirsky(1971)]{Mirsky}
L.~Mirsky.
\newblock A dual of {D}ilworth's decomposition theorem.
\newblock \emph{Amer. Math. Month.}, 78\penalty0 (8):\penalty0 876--877, 1971.
\newblock \doi{10.2307/2316481}.

\bibitem[Peles(1963)]{Peles}
M.~A. Peles.
\newblock On {D}ilworth's theorem in the infinite case.
\newblock \emph{Israel J.~Math.}, 1\penalty0 (2):\penalty0 108--109, 1963.
\newblock \doi{10.1007/BF02759806}.

\bibitem[Pouzet(1979)]{pouzet79}
M.~Pouzet.
\newblock Sur les cha\^ines d'un ensemble partiellement bien ordonn\'e.
\newblock \emph{Publ. D\'ep. Math. Lyon}, 16\penalty0 (1):\penalty0 21--26,
  1979.
\newblock URL \url{http://www.numdam.org/article/PDML_1979__16_1_21_0.pdf}.

\bibitem[Rado(1954)]{rado54}
R.~Rado.
\newblock Partial well-ordering of sets of vectors.
\newblock \emph{Mathematika}, 1\penalty0 (2):\penalty0 89--95, 1954.
\newblock \doi{10.1112/S0025579300000565}.

\bibitem[Rathjen and Weiermann(1993)]{RaWe}
M.~Rathjen and A.~Weiermann.
\newblock Proof-theoretic investigations on {K}ruskal's theorem.
\newblock \emph{Ann. Pure App. Logic}, 60\penalty0 (1):\penalty0 49--88, 1993.
\newblock \doi{10.1016/0168-0072(93)90192-G}.

\bibitem[Schmerl(2002)]{Schmerl}
J.~H. Schmerl.
\newblock Obstacles to extending {M}irsky's theorem.
\newblock \emph{Order}, 19\penalty0 (2):\penalty0 209--211, 2002.
\newblock \doi{10.1023/A:1016541101728}.

\bibitem[Schmidt(1979)]{schmidt79}
D.~Schmidt.
\newblock \emph{Well-Partial Orderings and Their Maximal Order Types}.
\newblock Habilitationsschrift, Heidelberg, 1979.

\bibitem[Schmidt(1981)]{schmidt81}
D.~Schmidt.
\newblock The relation between the height of a well-founded partial ordering
  and the order types of its chains and antichains.
\newblock \emph{J.~Comb. Theory B}, 31\penalty0 (2):\penalty0 183--189, 1981.
\newblock ISSN 0095-8956.
\newblock \doi{10.1016/S0095-8956(81)80023-8}.

\bibitem[Schmitz and Schnoebelen(2011)]{SS-icalp11}
S.~Schmitz and {\relax Ph}.~Schnoebelen.
\newblock Multiply-recursive upper bounds with {Higman}'s lemma.
\newblock In \emph{Proc.\ ICALP 2011}, volume 6756 of \emph{Lect. Notes Comput.
  Sci.}, pages 441--452. Springer, 2011.
\newblock \doi{10.1007/978-3-642-22012-8_35}.

\bibitem[Van~der Meeren(2015)]{vandermeeren15}
J.~Van~der Meeren.
\newblock \emph{Connecting the Two Worlds: Well-partial-orders and Ordinal
  Notation Systems}.
\newblock PhD thesis, Universiteit Gent, 2015.

\bibitem[Van~der Meeren et~al.(2015)Van~der Meeren, Rathjen, and
  Weiermann]{VdMRaWe}
J.~Van~der Meeren, M.~Rathjen, and A.~Weiermann.
\newblock Well-partial-orderings and the big {V}eblen number.
\newblock \emph{Arch. Math. Logic}, 54\penalty0 (1--2):\penalty0 193--230,
  2015.
\newblock \doi{10.1007/s00153-014-0408-5}.

\bibitem[Weiermann(2009)]{weiermann2009}
A.~Weiermann.
\newblock A computation of the maximal order type of the term ordering on
  finite multisets.
\newblock In \emph{Proc.\ CiE 2009}, volume 5635 of \emph{Lect. Notes Comput.
  Sci.}, pages 488--498. Springer, 2009.
\newblock \doi{10.1007/978-3-642-03073-4_50}.

\bibitem[West(1982)]{west82}
D.~B. West.
\newblock Extremal problems in partially ordered sets.
\newblock In I.~Rival, editor, \emph{Ordered Sets}, volume~83 of \emph{NATO
  Advanced Study Institutes Series}, pages 473--521. Springer, 1982.
\newblock \doi{10.1007/978-94-009-7798-3_16}.

\bibitem[Wolk(1967)]{Wolk}
E.~S. Wolk.
\newblock Partially well ordered sets and partial ordinals.
\newblock \emph{Fund. Math.}, 60\penalty0 (2):\penalty0 175--186, 1967.

\end{thebibliography}

\end{document}